\documentclass{article}
\usepackage{diagrams}

\usepackage{amsmath,amssymb,amsthm,subfigure,epsfig}
\usepackage{graphicx}

\usepackage[all,2cell]{xy}

\UseAllTwocells

\title{Topological Deformation of Higher Dimensional Automata}

\author{Philippe Gaucher$^{(*)}$, Eric Goubault$^{(**)}$\\
${ }^{(*)}$Institut de Recherche Math\'ematique Avanc\'ee,
ULP et  CNRS, \\
7 rue Ren\'e Descartes,
67084 Strasbourg Cedex \\
${ }^{(**)}$LIST (CEA - Recherche Technologique) \\
DTSI-SLA, CEA F91191 Gif-sur-Yvette Cedex \\}

\newcommand{\mb}[1]{{\mathbf #1}}
\newcommand{\C}{\mathcal{C}}

\newcommand{\N}{\mathbb{N}}
\newcommand{\R}{\mathbb{R}}
\newcommand{\U}{\mathcal{U}}
\newcommand{\V}{\mathcal{V}}

\newcommand{\intr}[1]{\stackrel{\circ}{#1}}

\newcommand{\brm}[1]{\rm{\mathbf{#1}}}
\newcommand{\de}{\partial}
\newcommand{\p}\times
\renewcommand{\vec}{\overrightarrow}
\newcommand{\iso}{\cong}
\newcommand{\e}{\vec{e}}
\newcommand{\CW}{{\brm{CW}}}
\newcommand{\diCW}{{\brm{glCW}}}
\newcommand{\lpohaus}{{\brm{LPoHaus}}}

\newcommand{\lpohauspp}{{\brm{LPoHaus}}_{**}}
\newcommand{\haus}{{\brm{Haus}}}
\renewcommand{\top}{{\brm{Top}}}

\newcommand{\potoppp}{{\brm{PoTop}}_{**}}

\newcommand{\ei}{\underline{\iota}}
\newcommand{\ef}{\underline{\sigma}}

\newcommand{\Ho}{{\brm{Ho}}}
\newcommand{\thaus}{{\brm{HAUS}}}

\newcommand{\be}{\begin{equation}}
\newcommand{\ee}{\end{equation}}
\newcommand{\bea}{\begin{eqnarray}}
\newcommand{\eea}{\end{eqnarray}}
\newcommand{\beas}{\begin{eqnarray*}}
\newcommand{\eeas}{\end{eqnarray*}}

\newtheorem{thm}{Theorem}[section]
\newtheorem{prop}[thm]{Proposition}
\newtheorem{lem}[thm]{Lemma}
\newtheorem{conj}[thm]{Conjecture}
\newtheorem{ex}[thm]{Example}

\newtheorem{question}[thm]{Question}
\newtheorem{cor}[thm]{Corollary}
\newtheorem{rem}[thm]{Remarque}

\newtheorem{defn}[thm]{Definition}
\newtheorem{propdef}[thm]{Proposition and Definition}

\newcommand{\bd}{\begin{defn}}
\newcommand{\ed}{\end{defn}}
\newcommand{\bcd}{\begin{defn}}
\newcommand{\ecd}{\end{defn}}

\newcommand{\bp}{\begin{prop}}
\newcommand{\ep}{\end{prop}}
\newcommand{\bth}{\begin{thm}}
\renewcommand{\eth}{\end{thm}}
\newcommand{\br}{\begin{rem}}
\newcommand{\er}{\end{rem}}
\newcommand{\bpf}{\begin{proof}}
\newcommand{\epf}{\end{proof}}

\newcommand{\fl}[1]{\ar@{->}[l]_{#1}}
\newcommand{\fr}[1]{\ar@{->}[r]^{#1}}
\newcommand{\fd}[1]{\ar@{->}[d]_{#1}}
\newcommand{\fu}[1]{\ar@{->}[u]^{#1}}
\newcommand{\f}[2]{\ar@{->}[#1]|{#2}}

\newcommand{\ff}[2]{\ar@2{->}[#1]|{#2}}

\newcommand{\vI}{\vec{I}}

\newcommand{\I}{\mathbb{I}}
\renewcommand{\P}{\mathbb{P}}

\renewcommand{\leq}{\leqslant}
\renewcommand{\geq}{\geqslant}

\def\limind{\setbox1=\hbox{\oalign{\vadjust{\vskip -2pt}%
     \rm lim\cr \vadjust{\vskip -2pt}
      \hidewidth$-\mkern -12mu\rightarrow$\hidewidth\cr}}
       \mathop{\box1}}

\date{July 2002}

\begin{document}

%\classification{55P15,55U05,68Q85}

%\keywords{homology, homotopy, concurrency, cubical set, CW-complex, higher dimensional automata, category, localization, partial order, partially ordered space}

\maketitle

\begin{abstract}
A local po-space is a gluing of topological spaces which are equipped with a
closed partial ordering representing the time flow. They are used as a
formalization of higher
dimensional automata (see for instance \cite{LFEGMRAlgebraic})
which model concurrent systems in computer science.
It is known \cite{ConcuToAlgTopo} that there are two distinct notions
of deformation of higher dimensional automata, ``spatial'' and ``temporal'',
leaving invariant
computer scientific properties like presence or absence of deadlocks. Unfortunately,
the
formalization of these notions is still unknown in the general case of
local po-spaces.

We introduce here a particular kind of local po-space, the
``globular CW-complexes'', for which we formalize these notions
of deformations and which are sufficient to formalize higher
dimensional automata. The existence of the category of globular
CW-complexes was already conjectured in \cite{ConcuToAlgTopo}.

After localizing the category of globular CW-complexes by
spatial and temporal deformations, we get a category (the category
of dihomotopy types) whose objects up to isomorphism represent
exactly the higher dimensional automata up to deformation. Thus
globular CW-complexes provide a rigorous mathematical
foundation to study from an algebraic topology point of view higher
dimensional automata and concurrent computations.
\end{abstract}

\tableofcontents

\section{Introduction}

\label{introduction}

Algebraic topological models have been used now for some years in concurrency
theory (concurrent database systems and fault-tolerant distributed systems as
well). The earlier models,
{\em progress graph} (see \cite{CoElSh71} for instance) have actually appeared in
operating systems theory, in particular for describing the problem
of ``deadly embrace''\footnote{as E. W. Dijkstra originally put it in
\cite{EWDCooperating}, now more usually called deadlock.} in
``multiprogramming systems''.

The basic idea is to give a description of what can happen when
several processes are modifying shared resources. Given a shared
resource $a$, we see it as its associated semaphore that rules
its behaviour with respect to processes. For instance, if $a$ is
an ordinary shared variable, it is customary to use its semaphore
to ensure that only one process at a time can write on it (this
is mutual exclusion). A semaphore is nothing but a register which counts
the number of times a shared object can still be accessed by processes. In
the case of usual shared variables, this register is initialized with value 1,
processes trying to access (read or write) on the corresponding variable 
compete in order to get it first, then the semaphore value is decreased:
we say that the semaphore has been locked\footnote{Of course this operation
must be done ``atomically'', meaning that the semaphore itself must be handled
in a mutually exclusive manner: this is done at the hardware level.} by the process. 
When it is equal to zero, all processes trying to access this semaphore are blocked, waiting
for the process which holds the lock to relinquish it, typically when it
has finished reading or writing on the corresponding variable: 
the value of the semaphore is then increased. 

When the semaphores are initialized
with value one, meaning that they are associated with shared variables accessed in
a mutually exclusive manner, they are called binary semaphores. When a shared data (identified with its semaphore) can be accessed by one or more processes, meaning
that the corresponding semaphore has been initialized with a value greater than one,
it is called a counting semaphore. 

Given $n$ deterministic sequential
processes $Q_1,\ldots,Q_n$, abstracted as a sequence of locks and
unlocks on (semaphores associated with) shared objects, $ Q_i=R^1 a_i^1.R^2 a_i^2 \cdots
R^{n_i} a_i^{n_i} $ ($R^k$ being $P$ or $V$\footnote{Using E. W.
Dijkstra's notation $P$ and $V$ \cite{EWDCooperating} for
respectively acquiring and releasing a lock on a semaphore.}),
there is a natural way to understand the possible behaviours of
their concurrent execution, by associating to each process a
coordinate line in $\R^n$. The state of the system corresponds to
a point in $\R^n$, whose $i$th coordinate describes the state (or
``local time'') of the $i$th processor.

\begin{figure}
\begin{center}
\includegraphics[width=7cm]{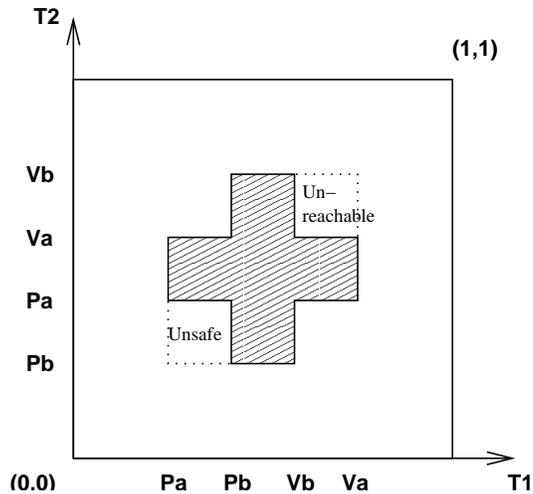}
\end{center}
\caption{Example of a progress graph}
\label{progress1}
\end{figure}

Consider a system with finitely many processes running altogether.
We assume that each process starts at (local
time) 0 and finishes at (local time) 1; the $P$ and $V$ actions correspond to
sequences of real numbers between 0 and 1, which reflect the order of
the $P$'s and $V$'s. The initial state is $(0,\dots ,0)$ and the final
state is $(1,\dots ,1)$.
An example consisting of the two processes
$T_1=P a.P b.V b.Va$ and $T_2=P b.P a.V a.V b$ gives rise to the two
dimensional {\em progress graph\/} of Figure \ref{progress1}.

The shaded area represents states which are not allowed in any execution
path, since they correspond to mutual exclusion. Such states constitute the
{\em forbidden area}. An {\em execution path} is a path from the initial
state $(0,\ldots,0)$ to the final state $(1,\ldots,1)$ avoiding the forbidden area and
increasing in each coordinate - time cannot run backwards. We call
these paths {\em directed paths} or dipaths. This entails that
paths reaching the
states in the dashed square underneath the forbidden region, marked ``unsafe''
are deemed to deadlock, i.e. they cannot possibly reach the allowed terminal
state which is $(1,1)$ here. Similarly, by reversing the direction of time,
the states in the square above the forbidden region, marked ``unreachable'',
cannot be reached from the initial state, which is $(0,0)$ here.
Also notice that all terminating paths above the
forbidden region are ``equivalent'' in some sense, given that they are all
characterized by the fact that $T_2$ gets $a$ and $b$ before $T_1$ (as far as
resources are concerned, we call this a {\em schedule}). Similarly, all paths
below the forbidden region are characterized by the fact that $T_1$ gets
$a$ and $b$ before $T_2$ does.

\begin{figure}
\begin{center}
\includegraphics[width=7cm]{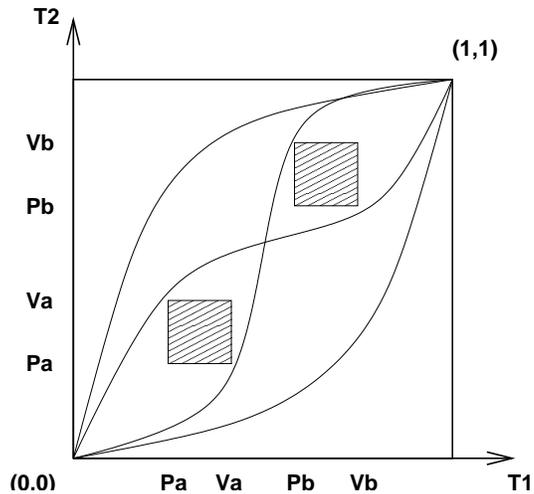}
\end{center}
\caption{The progress graph corresponding to $Pa. Va. Pb. Vb \mid
Pa. Va. Pb. Vb$}
\label{shape1}
\end{figure}

\begin{figure}
\begin{center}
\includegraphics[width=7cm]{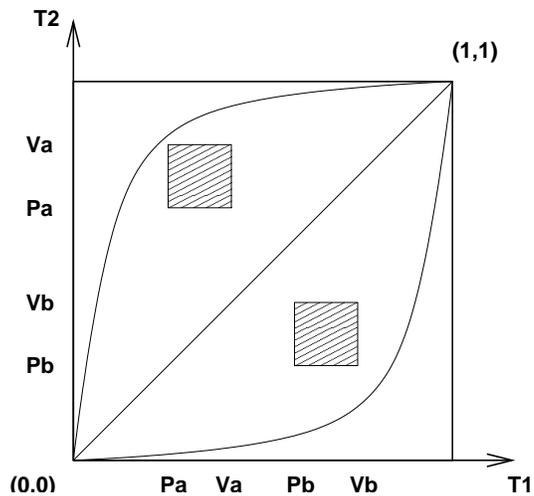}
\end{center}
\caption{The progress graph corresponding to $Pb. Vb. Pa. Va \mid
Pa. Va. Pb. Vb$}
\label{shape2}
\end{figure}

On this picture, one can already recognize many ingredients that are at the center
of the main problem of algebraic topology, namely the classification of shapes
modulo ``elastic deformation''. As a matter of fact, the actual coordinates
that are chosen for representing the times at which $P$s and $V$s occur
are unimportant, and these can be ``stretched'' in any manner, so the properties
(deadlocks, schedules etc.) are invariant under some notion of deformation,
or {\em homotopy}. This is only a particular kind of homotopy though, and this
explains why a new theory has to be designed. We call it (in
subsequent work) {\em directed homotopy} or {\em dihomotopy}
in the sense that it should preserve
the direction of time. For instance, the two homotopic shapes, all of which
have two holes, of Figure \ref{shape1} and Figure \ref{shape2} have a
different number of dihomotopy classes of dipaths. In Figure \ref{shape1} there
are essentially four dipaths up to dihomotopy (i.e. four schedules corresponding
to all possibilities of accesses of resources $a$ and $b$) whereas in Figure
\ref{shape2}, there are essentially three dipaths up to dihomotopy.

Progress graphs have actually a nice topological model; they
are compact order-Hausdorff spaces (see \cite{nachbin},
\cite{johnstone}), i.e.  are compact Hausdorff spaces with a closed
(global) partial order.  More general topological models are needed in
general, in which the partial order is only defined locally, and have
been introduced and motivated in \cite{LFMRDetecting},
\cite{LFEGMRDetecting} and \cite{LFEGMRAlgebraic}.  The precise
definitions and properties are given in Section~\ref{pospace}.

The natural
combinatorial
notion which discretizes this topological framework
is that of a {\em precubical set},
which is a collection of points (states), edges (transitions), squares, cubes
and hypercubes (higher-dimensional transitions representing the truly-concurrent
execution of some number of actions). This is introduced in \cite{VPModeling}
as well as possible formalizations using $n$-categories, and a notion
of homotopy. 
These precubical sets are called
{\em Higher-Dimensional Automata} (HDA) following \cite{VPModeling} because
it really makes sense to consider a hypercube as some form of transition
(as in transition systems, used in semantics of programming languages).
We show the precise relation between this model and the new topological
model we introduce here (``globular CW-complexes'') in Section~\ref{globular},
the relation between local po-spaces and cubical sets can be found in
\cite{LFEGMRAlgebraic}.

There are other formulations of the same problems using homological
methods \cite{HDA}, strict globular $\omega$-categories \cite{Gau}. An
important motivation in these pieces of work is that of ``reducing the
complexity'' of the semantics (given by a local po-space for instance)
by considering deformation retracts. The classification of the
possible concurrent semantics (and behaviours) should then be the
result finding the right ``(di-)homotopy types''. This calls for a
suitable notion of (di-)homotopy equivalence, and for starting with a
reasonable category of local po-spaces.  In the case of ordinary
homotopy theory, we have to restrict to the category of CW-complexes;
the category of topological spaces being far too big for practical
purposes. The situation is even worse here, we not only have to
restrict on the topology part, but also on the local po-structures.

We give in this paper a notion of CW-complex, called globular
CW-complex which meets the basic requirements of what we expect to be
a ``directed cellular complex''. It has been obtained by
mimicking the well-known concept of CW-complexes, but built from
``directed'' cells.  This is the purpose of
Section~\ref{section_diCW}. Still in the same section, we introduce
the fundamental functor called the \textit{globe functor}, from the category of
topological spaces to the category of po-spaces. This functor is the
key to understanding how things work in the directed situation. In
particular, it yields an embedding of the category of homotopy types
into the new category of dihomotopy types (Theorem~\ref{embed}). This
embedding has a lot of important consequences that are sketched in the
perspectives section of \cite{Ditype}.

Once the right notion has been given, we make explicit the link
between the globular CW-complexes and some geometric notions above
mentioned, that is the local po-spaces and the precubical sets in
Section~\ref{lien}. We prove that every globular CW-complex
is a local po-space indeed (Theorem~\ref{CWpo}) and that there exists
a geometric realization of any precubical set as a globular CW-complex
(subsection~\ref{appA}).

Next in Section~\ref{ST} we recast in the globular CW-complex
framework the notion of spatial and temporal deformations informally
presented in \cite{ConcuToAlgTopo} and whose consequences are
informally explored in \cite{Ditype}. For that we construct, by
localization of the category of globular CW-complexes with respect to
appropriate morphisms, three categories whose isomorphism classes of
objects are exactly the globular CW-complexes modulo spatial
deformations (Theorem~\ref{HOS}), the globular CW-complexes modulo
temporal deformations (Theorem~\ref{HOT}) and at last the globular
CW-complexes modulo spatial and temporal deformations together
(Theorem~\ref{HOST}). This latter category will be called the category
of \textit{dihomotopy types}.

Then  Section~\ref{embed-ho-diho} is devoted to making explicit
the link between homotopy types and dihomotopy types. The
introduction of the \textit{path space} of a globular CW-complex
between two points of its $0$-skeleton is the essential
ingredient in the proof of Theorem~\ref{th} and
Corollary~\ref{bij-hom-dihom}. This allows us to derive the
embedding theorem Theorem~\ref{embed} which states that there
exists an embedding of the category of homotopy types into that
of dihomotopy types.  This notion of path spaces also allows us
to provide a conjectural statement for the analogue of the
Whitehead theorem in the directed framework, and to check it in
the case of globes.

Section~\ref{why} focuses on a very striking reason why it is necessary to
work with ``non-contrac\-ting'' maps everywhere. It was not really possible
to justify this axiom while the definition of globular CW-complexes
was being given in Section~\ref{section_diCW}. Only one reason is
described. Indeed there are lots of other algebraic reasons which are
out of the scope of this paper.

\section{Globular CW-complexes}\label{section_diCW}

\label{globular1}

This section is devoted to the introduction of the category $\diCW$ of
globular CW-complexes and to the comparison with the usual notion of
CW-complex.

\subsection{Closed partial order}

\bd Let $X$ be a topological space. A binary relation $R$ on $X$ is
closed if the graph of $R$ is a closed subset of the cartesian product
$X\p X$. \ed

It is a well-known fact that any topological space $X$ endowed with a
closed partial order is necessarily Hausdorff (see for instance
\cite{nachbin}, \cite{johnstone}).

\bd A pair $(X,\leq_X)$ where $X$ is a topological space and $\leq_X$ a closed
partial order is called a \textit{global po-space}.  \ed

In most cases, the partial order of a global po-space $X$ will
be simply denoted by $\leq$. Here are two fundamental examples of
global po-spaces for the sequel\thinspace:
\begin{enumerate}
\item The \textit{achronal segment} $\I$ is defined to be the segment
$[0,1]$ endowed with the closed partial ordering $x\leq_{\I} y$ if and
only if $x=y$.
\item The \textit{directed segment} $\vI$ is defined to be the segment
$[0,1]$ endowed with the closed partial ordering $x\leq_{\vI} y$ if
and only if $y-x$ is non-negative.
\end{enumerate}

We will describe gluings of global po-spaces (i.e. local po-spaces)
in
Section~\ref{pospace}.

\subsection{The globe of a topological space}

Let $n\geq 1$. Let $D^n$ be the closed $n$-dimensional disk defined by
the set of points $(x_1,\dots,x_n)$ of $\R^n$ such that $x_1^2+\dots
+x_n^2\leq 1$ endowed with the topology induced by that of $\R^n$. Let
$S^{n-1}=\de D^n$ be the boundary of $D^n$ for $n\geq 1$, that is the
set of $(x_1,\dots,x_n)\in D^n$ such that $x_1^2+\dots +x_n^2=1$.
Notice that $S^0$ is the discrete two-point topological space
$\{-1,+1\}$.  Let $D^0$ be the one-point topological
space. And let $e^n:=D^n-S^n$.

\begin{figure}
\begin{center}
\includegraphics[width=7cm]{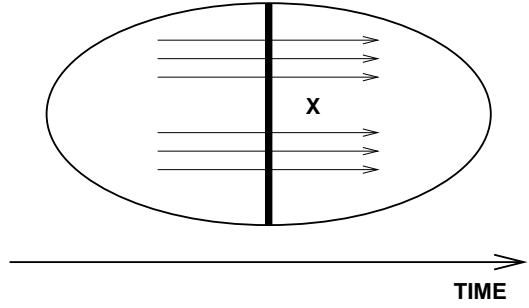}
\end{center}
\caption{Symbolic representation of $Glob(X)$ for some topological space $X$}
\label{exglob}
\end{figure}

The fundamental ingredient in all further constructions is the
globe functor (Figure~\ref{exglob}) defined below, which will give
rise to a particular family of global po-spaces. Loosely speaking
the globe of a topological space $X$ is the reduced suspension of $X$
equipped with some closed partial ordering representing the time flow.

The underlying topological space of the \textit{globe} of a
topological space $X$, $Glob(X)$, is therefore the quotient of the
product space $X\p [0,1]$ by the relations $(x,0)=(x',0)$ and
$(x,1)=(x',1)$ for any $x,x'\in X$. By convention, the equivalence
class of $(x,0)$ (resp.  $(x,1)$) in $Glob(X)$ will be denoted by
$\ei$ (resp. $\ef$). We can partially order $Glob(X)$ using the
standard order $\leq_I$ on $I=[0,1]$ as follows\thinspace:

\bp\label{glob} Let $X$ be a Hausdorff topological space and consider the
partial ordering of $X\p I$ defined by
\(\mathcal{R}=\{((x,t),(x,t')), (x,t,t')\in X\p I\p I\hbox{ and }t\leq_I t'\}\).
Then its image by the canonical
surjection $s$ from $X\p I$ to $Glob(X)$ is a closed partial ordering on
$Glob(X)$. \ep

The partial order relation on $Glob(X)$ is as follows:
\begin{itemize}
\item $(x,0) \leq (x',t')$ for all $x, x', t' \in X \times X \times
I$,
\item when $t, t' \in ]0,1[ \times ]0,1[$, $(x,t) \leq (x',t')$
if and only if $x=x'$,
\item $(x',t') \leq (x,1)$ for all $x, x', t' \in X \times X \times
I$.
\end{itemize}

\bpf By the
homeomorphism $(x,t,x',t')\mapsto (x,x',t,t')$ from $X\p I \p X \p I$
to $X\p X \p I \p I$, one sees that $\mathcal{R}$ is a closed subset
of $X\p I \p X \p I$ if and only if $Diag(X)\p \{(t,t')\in I\p I,t\leq
t'\}$ is a closed subset of $X\p X \p I \p I$ where $Diag(X)$ is the
diagonal $\{(x,x) / x \in X\}$ of $X$.
Since $X$ is Hausdorff, then its diagonal is closed
and $\mathcal{R}$ as well. By definition of the quotient topology,
$s(\mathcal{R})$ is closed if and only if $s^{-1}\circ s(\mathcal{R})$
is a closed subset of $X\p I$.  It suffices then to notice that
\(s^{-1}\circ s(\mathcal{R})=\left((X\p\{0\})\p (X\p I)\right) \cup
\left((X\p I)\p (X\p \{1\})\right)\cup \mathcal{R}\) to complete the
proof.  \epf

Ordinary CW-complexes are built by gluing cells $e^n$. We want to
define a {\em directed version} of CW-complexes, and the simplest way
to do so is to add a ``direction'' to these cells. Every HDA can be
seen indeed as a pasting of $n$-cubes or of $n$-globes, depending on whether
one chooses the cubical approach or the globular approach to model the
execution paths, the higher dimensional homotopies between them and
the compositions between them (see Section
\ref{introduction} and \cite{EGUsers} for more explanations). Loosely 
speaking, directed cells such as $Glob(e^n)$ are ``equivalent'' modulo 
directed deformations to a $n$-cube with the usual cartesian partial 
ordering defined by $(x_1,\dots,x_n)\leq (y_1,\dots,y_n)$ if and only 
if for any $i$, $x_i\leq y_i$.

\subsection{Globular CW-complex\thinspace: definition and examples}

Let $\vec{D}^{n+1}:=Glob(D^n)$ and $\vec{S}^{n+1}:=Glob(S^n)$ for
$n\geq 0$.  Notice that there is a canonical inclusion of global
po-spaces $\vec{S}^{n}\subset \vec{D}^{n+1}$ for $n\geq 1$. By
convention, let $\vec{S}^0:=\{0,1\}$ with the trivial ordering
($0$ and $1$ are not comparable). There is a canonical inclusion
$\vec{S}^0\subset \vec{D}^{1}$ such that $0$ is mapped onto $\ei$ 
(the initial state of $\vec{D}^1$) and
$1$ is mapped onto $\ef$ (the final 
state of $\vec{D}^1$), which is a morphism of
po-spaces.

\bd\label{ndicell} For any $n\geq 1$, $\vec{D}^{n}-\vec{S}^{n-1}$
together with the closed partial ordering induced by $I$ is called the
$n$-dimensional globular cell. More generally, every pair $(X,\leq)$,
where $X$ is a topological space and $\leq$ a closed partial ordering on
$X$, isomorphic to  $\vec{D}^{n}-\vec{S}^{n-1}$ for some $n$ will be called a
$n$-dimensional globular cell. \ed

Now we are going to describe the process of attaching globular cells.

\begin{enumerate}
\item Start with a discrete set of points $X^0$.
\item Inductively, form the $n$-skeleton $X^n$ from $X^{n-1}$
by attaching globular $n$-cells $\vec{e}^n_\alpha$ via maps
$\phi_\alpha:\vec{S}^{n-1}\longrightarrow X^{n-1}$ with
$\phi_\alpha(\ei),\phi_\alpha(\ef)\in X^0$ such that\footnote{This
condition will appear to be necessary in the sequel.}\thinspace: for every
non-decreasing map $\phi$ from $\vI$ to $\vec{S}^{n-1}$ such that
$\phi(0)=\ei$ and $\phi(1)=\ef$, there exists $0=t_0<\dots
<t_k=1$ such that $\phi_\alpha\circ \phi(t_i)\in X^0$ for any
$0\leq i\leq k$ which must satisfy
\begin{enumerate}
\item for any $0\leq i\leq k-1$, there exists a globular cell of
dimension $d_i$ with $d_i\leq n-1$ $\psi_i:\vec{D}^{d_i}\rightarrow
X^{n-1}$ such that for any $t\in[t_i,t_{i+1}]$, $\phi_\alpha\circ
\phi(t)\in \psi_i(\vec{D}^{d_i})$\thinspace;
\item for $0\leq i\leq k-1$, the restriction of $\phi_\alpha\circ
\phi$ to $[t_i,t_{i+1}]$ is non-decreasing\thinspace;
\item the map $\phi_\alpha\circ \phi$ is non-constant\thinspace;
\end{enumerate}
Then $X^n$ is the quotient space of the disjoint union
$X^{n-1}\bigsqcup_{\alpha}\vec{D}^n_\alpha$ of $X^{n-1}$ with a
collection of $\vec{D}^n_\alpha$ under the identification $x\sim
\phi_\alpha(x)$ for $x\in \vec{S}^{n-1}_\alpha\subset
\de\vec{D}^n_\alpha$. Thus as set,
$X^n=X^{n-1}\bigsqcup_{\alpha}\vec{e}^n_\alpha$ where each
$\vec{e}^n_\alpha$ is a $n$-dimensional globular cell.
\item One can either stop this inductive process at a finite stage,
  setting $X=X^n$, or one can continue indefinitely, setting
  $X=\bigcup_n X^n$. In the latter case, $X$ is given the weak
  topology\thinspace: a set $A\subset X$ is open (or closed) if and only if
  $A\cap X^n$ is open (or closed) in $X^n$ for some $n$ (this topology
  is nothing else but the direct limit of the topology of the $X^n$,
$n \in \N$). Such a $X$ is called a
  globular CW-complex and $X_0$ and the collection of
  $\vec{e}^n_\alpha$ and its attaching maps
  $\phi_\alpha:\vec{S}^{n-1}\longrightarrow X^{n-1}$ is called the
  cellular decomposition of $X$.
\end{enumerate}

As for usual CW-complexes (see \cite{Hat} Proposition A.2.),
a globular cellular decomposition of a given
globular CW-complex $X$ yields characteristic maps
$\phi_\alpha:\vec{D}^{n_\alpha}\rightarrow X$ satisfying\thinspace:
\begin{enumerate}
\item The mapping
$\phi_\alpha\restriction_{\vec{D}^{n_\alpha}-\vec{S}^{n_\alpha-1}}$
induces an homeomorphism from $\vec{e}^{n_\alpha}$ to its image.
\item All the previous globular cells are disjoint and their union
gives back $X$.
\item A subset of $X$ is closed if and only if it meets the closure
of each globular cells of $X$ in a closed set.
\end{enumerate}

We will consider without further mentioning that the segment $\vI$ is
a globular CW-complex, with $\{0,1\}$ as its $0$-skeleton.

\begin{propdef}\label{def_exec_path}
Let $X$ be a globular CW-complex with characteristic maps
$(\phi_\alpha)$. Let $\gamma$ be a continuous
map from $\vI$ to $X$. Then $\gamma([0,1])\cap X^0$ is finite.
Suppose that there exists $0\leq t_0<\dots<t_n\leq 1$ with
$n\geq 1$ such that $t_0=0$, $t_n=1$,  such that
for any $0\leq i\leq n$, $\gamma(t_i)\in X^0$, and at last
such that for any $0\leq i\leq n-1$,
there exists an $\alpha_i$ (necessarily unique) such that for $t\in [t_i,t_{i+1}]$,
$\gamma(t)\in \phi_{\alpha_i}(\vec{D}^{n_\alpha})$. Then such a
$\gamma$ is called an \textit{execution path} if
the restriction $\gamma\restriction_{[t_i,t_{i+1}]}$
is non-decreasing.
\end{propdef}

\bpf Obvious. \epf

By constant execution paths, one means an execution paths $\gamma$ such that
$\gamma([0,1])=\{\gamma(0)\}$ with $\gamma(0) \in X^0$. The points (i.e. elements of the $0$-skeleton)
of a given globular CW-complexes $X$ are also called
\textit{states}. Some of them are fairly special:

\bd Let $X$ be a globular CW-complex. A point $\alpha$ of $X^0$ is
initial (resp. final) if for any execution path $\phi$ such that
$\phi(1)=\alpha$ (resp. $\phi(0)=\alpha$), then $\phi$ is the constant
path $\alpha$. \ed

\bp\label{globe-CW} If $X$ is a CW-complex, then $Glob(X)$ is a globular
CW-complex by setting
\[Glob(X)^0=\{\ei,\ef\}\] for $x\in X$.
\ep

\bpf Since $X$ is a CW-complex, then it is described by cells and
attaching maps. There exists topological spaces $X^n$ with
$X=\bigcup_n X^n$ with the weak topology and
$\phi_\alpha:{S}^{n-1}\longrightarrow X^{n-1}$ (for $\alpha$
belonging to some set of indexes) continuous maps which describe
how to go from $X^{n-1}$ to $X^n$; we have the following
co-cartesian diagram of topological spaces:
$$
\begin{diagram}
{S}^{n-1} & \rTo^{\phi_\alpha} & X^{n-1} \\
\dTo^{i} & & \dTo \\
D^{n} & \rTo & X^n \\
\end{diagram}
$$
where $i$ is the inclusion of $S^{n-1}$ into $D^n$ as its boundary $\de{D}^n$.

Let us describe inductively $Glob(X)$ as a globular CW-complex. We
begin by setting $Glob(X)^0=\{\ei,\ef\}$.
Then we apply inductively the functor $Glob(-)$ on the co-cartesian diagram above, to obtain
a new co-cartesian diagram by Theorem~\ref{dlim}\thinspace:
$$
\begin{diagram}
Glob({S}^{n-1})\iso \vec{S}^n & \rTo^{Glob(\phi_\alpha)} & Glob(X^{n-1}) \\
\dTo^{Glob(i)} & & \dTo \\
Glob(D^{n})\iso \vec{D}^{n+1} & \rTo & Glob(X^n) \\
\end{diagram}
$$
First of all, it is easy to see that $Glob(i)$ induces a homeomorphism
from $\vec{S}^n$ onto the boundary $\de{\vec{D}}^{n+1}$ of $\vec{D}^{n+1}$,
therefore is the inclusion morphism we expect. We now have to check that
$Glob(\phi_{\alpha})$ is a correct attaching map for globular CW-complexes.
For $(x,u) \in Glob(S^{n-1})$ ($x \in S^{n-1}$, $u \in \vI$), we have
$Glob(\phi_{\alpha})(x,u)=(\phi_\alpha(x),u)$. We have to see that it is non-decreasing.
Let $(x,u)$ and $(x',u')$ be two elements of $Glob(S^{n-1})$ such that
$(x,u) \leq (x',u')$. We have the following cases:
\begin{itemize}
\item $u=0$ then $Glob(\phi_{\alpha})(x,u)=\ei$, thus is less or equal than
$Glob(\phi_{\alpha})(x',u')$,
\item $u'=1$ then $Glob(\phi_{\alpha})(x',u')=\ef$, thus is greater or equal than
$Glob(\phi_{\alpha})(x,u)$,
\item $0 < u < 1$ (the case $u=1$ is trivial since it implies $u'=1$, which is
the previous case)
then $u \leq u'$ and $x=x'$. Thus,
$Glob(\phi_{\alpha})(x,u)=(\phi_{\alpha}(x),u) \leq (\phi_{\alpha}(x'),u)=
Glob(\phi_{\alpha})(x',u')$.
\end{itemize}

That $Glob(\phi_\alpha)$ is non-contracting is due to the fact that
$Glob(\phi_\alpha)(\ei)\neq Glob(\phi_\alpha)(\ef)$.

\epf

\bp Every globular CW-complex is a CW-complex.  \ep

\bpf This is due to the fact that $\vec{e}^n_\alpha$ is homeomorphic to
${e}^n_\alpha$. \epf

We end up the section with the Swiss Flag example of Figure~\ref{progress1} 
described as a globular CW-complex. We will use this example throughout
the article to illustrate different constructions and ideas.

\begin{figure}
\begin{center}
\includegraphics[width=7cm]{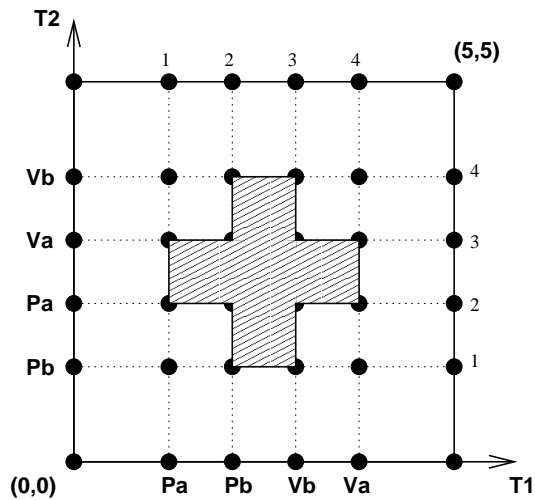}
\end{center}
\caption{Running example of globular CW-complex}
\label{progress10}
\end{figure}

\begin{ex}\label{runningexample}
Consider the discrete set $SW^0=\{0,1,2,3,4,5\}\p \{0,1,2,3,4,5\}$.
Let 
\beas 
\mathcal{S}&=&\left\{((i,j),(i+1,j))\hbox{ for } (i,j)\in\{0,\dots,4\}\p \{0,\dots,5\}\right\}\\
&\cup& \left\{((i,j),(i,j+1))\hbox{ for } (i,j)\in\{0,\dots,5\}\p \{0,\dots,4\}\right\}\\
&\backslash & \left(
\{((2,2),(2,3)),((2,2),(3,2)), ((2,3),(3,3)),((3,2),(3,3))\}
\right)
\eeas
The globular CW-complex $SW^1$ is obtained from $SW^0$ by attaching 
a copy of $\vec{D}^1$ to each pair $(x,y)\in \mathcal{S}$ with 
$x\in SW^0$ identified with $\ei$ and $y\in SW^0$ identified with $\ef$. 
The globular CW-complex $SW^2$ is obtained from $SW^1$ by attaching 
to each square $((i,j),(i+1,j+1))$ except $(i,j)\in\{(2,1),(1,2),(2,2),(3,2),(2,3)\}$ a globular cell $\vec{D}^2$ such that 
each execution path $((i,j),(i+1,j),(i+1,j+1))$ and 
$((i,j),(i,j+1),(i+1,j+1))$ is identified with one of the execution 
path of $\vec{S}^1$ (there is no unique choice to do that). Let $SW=SW^2$ 
(cf. Figure~\ref{progress10} where the bold dots represent the points 
of the $0$-skeleton). 
\end{ex}

\subsection{Morphism of globular CW-complexes}
\label{globCWcat}

\bd\label{morphism-glCW} The category $\diCW$ of globular CW-complexes
is the category having as objects the globular CW-complexes and as
morphisms the continuous maps  $f:X\longrightarrow Y$
satisfying the
two following properties\thinspace:
\begin{itemize}
\item $f(X^0)\subset Y^0$
\item for every non-constant execution path $\phi$ of $X$,
$f\circ \phi$ must not only be an execution path ($f$ must preserve
partial order), but also $f\circ \phi$ must be non-constant
as well \thinspace: we say that $f$ must be \textit{non-contracting}.
\end{itemize}
\ed

The condition of non-contractibility is very analogous to the notion
of non-contracting $\omega$-functors appearing in \cite{Gau}.  Notice
also that the attaching maps in the definition of a globular
CW-complex are morphisms in $\diCW$. This non-contractibility
condition will be justified in Section~\ref{why}.

A non-constant execution path of a globular CW-complex $X$ induces a
morphism of globular CW-complexes from $\vI$ to $X$.

\bp\label{globe-CW2} The functor $Glob(-)$ induces a
functor still denoted by $Glob(-)$ from the category $\CW$ of
CW-complexes to the category $\diCW$ of globular CW-complexes.
\ep

\bpf  It is an immediate consequence of Proposition~\ref{globe-CW}. \epf

\begin{figure}
\begin{center}
\includegraphics[width=7cm]{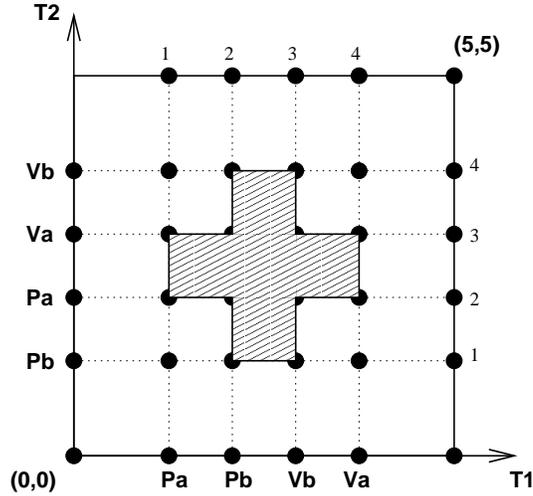}
\end{center}
\caption{Another example of globular CW-complex}
\label{progress11}
\end{figure}

\begin{ex}\label{ex_T}
As example of morphisms of globular CW-complexes, the identity 
map of $[0,5]\p [0,5]$ induces a morphism of globular CW-complexes 
from the globular CW-complex of Figure~\ref{progress11} 
to that of Figure~\ref{progress10}. It will be an example of 
T-homotopy equivalence (Section~\ref{Tequiv}).  \end{ex}

\section{Relation with other formalizations}
\label{lien}

\subsection{Gluing closed partial orderings}
\label{pospace}

Let us remind some definitions to fix the notations. The category of
Hausdorff topological spaces with the continuous maps as morphisms
will be denoted by $\haus$. The category of general topological spaces
without further assumption will be denoted by $\top$.

Let $(X,R)$ be a global po-space.  We say that $(U,\leq)$ is a
sub-po-space of $(X,R)$ if and only if it is a po-space such that $U$
is a topological subspace of $X$ and such that $\leq$ is the
restriction of $R$ to $U$.

Let $X$ be a Hausdorff topological space.
A collection $\U(X)$ of po-spaces $(U,\leq_U)$ covering $X$ is called a
\textit{local po-structure} if for every $x\in X$, there exists a
po-space $(W(x),\leq_{W(x)})$ such that:
\begin{itemize}
\item
 $W(x)$ is an open neighborhood
containing $x$,
\item the restrictions of $\leq_U$ and
$\leq_{W(x)}$ to $W(x)\cap U$ coincide for all $U\in \U(X)$ such
that $x\in U$. This can be stated as: $y\leq_U z$ if and only if
$y\leq_{W(x)} z$ for all $U\in \U(X)$ such that $x\in U$ and for
all $y,z\in W(x)\cap U$.  Sometimes, $W(x)$ will be denoted by
$W_X(x)$ to avoid ambiguities. Such a $W_X(x)$ is called a
po-neighborhood.
\end{itemize}

Two local po-structures are equivalent if their union is a local
po-structure. This defines an equivalence relation on the set of local
partial structures of $X$. A topological space together with an
equivalence class of local po-structures is called a \textit{local
  po-space} \cite{LFEGMRAlgebraic}. Notice that a global po-space is a
local po-space.

A morphism $f$ of local po-spaces (or \textit{dimap}) from
$(X,\U)$ to $(Y,\V)$ is a continuous map from $X$ to $Y$ such
that for every $x\in X$,
\begin{itemize}
\item
there is a po-neighborhood $W(f(x))$ of $f(x)$ in $Y$,
\item
there exists a po-neighborhood $W(x)$ of $x$ in $X$ with
$W(x)\subset f^{-1}(W(f(x)))$,
\item for $y,z\in W(x)$, $y\leq z$ implies
$f(y)\leq f(z)$.
\end{itemize}

In particular, a dimap $f$ from a po-space $X$ to a
po-space $Y$ is a continuous map from $X$ to $Y$ such that for any
$y,z\in X$, $y\leq z$ implies $f(y)\leq f(z)$. A morphism $f$
of local po-spaces from $[0,1]$ endowed with the usual ordering
(denoted by $\vec{I}$)
to a local po-space $X$ is called \textit{dipath} or sometime
\textit{execution path}.

The category of Hausdorff local po-spaces with the dimaps as morphisms
will be denoted by $\lpohaus$. The mapping $Glob(-)$ of
Proposition~\ref{glob} yields a functor from $\haus$ to $\lpohaus$.

\subsection{Globular CW-complex and local po-space}

\label{globular}

We now relate globular CW-complexes with local po-spaces.

\paragraph{Convention}

In the
sequel, for any $X$ and $Y$ two topological spaces, we endow the disjoint
sum $X\sqcup Y$ with the final topology induced by both
inclusion maps $X\subset X\sqcup Y$ and $Y\subset X\sqcup Y$.

Both following lemmas summarize well-known facts about topological
spaces\thinspace: see \cite{Rotman} exercises 8.12 and 8.13.

\begin{lem}\label{ferm}
Let $\phi$ be a closed continuous map from $X$ to $Y$ and
let $Z\subset Y$. Let $U$ be an open subset of $X$ containing $\phi^{-1}(Z)$.
Then there exists an open subset $V$ of $Y$ such that $Z\subset V$ and
$\phi^{-1}(V)\subset U$.
\end{lem}

%\bpf Let $V:=Y-\phi(X-U)$. Since $\phi$ is closed, $V$ is a closed
%subset of $Y$. The inclusion $\phi^{-1}(V)\subset U$ is obvious. Now
%if $z\in Z$, then either $z\in Y-\phi(X)\subset V$, or $z=\phi(x)$ for
%some $x\in X$. And $x\in \phi^{-1}(Z)\subset U$ implies $x\in U$.
%Therefore $z\notin \phi(X-U)$.
%\epf

\begin{lem}\label{base}
Let $A$ be a closed subset of $X$. Let $f$ be a continuous map
from $A$ to $Y$. Consider the quotient topological space
$X \sqcup_f Y\thinspace:= (X\sqcup Y)/\sim$ where
$\sim$ is the equivalence relation on $X\sqcup Y$
generated by $\{(a,f(a))\in (X\sqcup Y)\p (X \sqcup Y),a\in A\}$.
Let $\phi$ be the canonical continuous map from $X\sqcup Y$ to
$X \sqcup_f Y$. Then
\begin{enumerate}
\item\label{deftop} A subset $\Omega$ of $X \sqcup_f Y$ is open (resp. closed) in $X \sqcup_f Y$ if and only if
$\phi^{-1}(\Omega)\cap X$ is open (resp. closed) in $X$ and $\phi^{-1}(\Omega)\cap Y$
is open (resp. closed) in $Y$.
\item\label{casferme} As soon as $A$ is a closed subset of $X$, $X-A$ can be identified to
the open subset $\phi(X-A)$ of $X \sqcup_f Y$ and $Y$ can be
identified to the closed subset $\phi(Y)$ of $X \sqcup_f Y$.
\item\label{cascompact} If $f(A)$ is a closed subset of $Y$, then $Y-f(A)$ can be
identified to the open subset $\phi(Y-f(A))$ of $X \sqcup_f Y$ and
$f(A)$ to the closed subset $\phi(f(A))$ of $Y$.
\item\label{cascompactencore} If $A$ is compact, then $\phi$ is a
closed map.
\item\label{separation} If $A$ is compact and if $X$ and $Y$ are Hausdorff, then
$X \sqcup_f Y$ is Hausdorff.
\end{enumerate}
\end{lem}

\bth\label{CWpo} Every globular CW-complex $X$ is a local po-space.
\eth

\bpf We prove that attaching globular $n$-cells to any locally compact local
po-space still defines a local po-space. As points are trivial local
po-spaces, the theorem will follow from an easy induction.

First we say that a local po-structure is \textit{small} if for all
$U$ and $V$ in the open covering defining the local po-structure,
$\leq_U$ and $\leq_V$ coincide on $U \cap V$. It is easy to see that
all local po-spaces $X$ admit (in its equivalence class of coverings)
a \textit{small local po-structure}\thinspace: if $W_X(x)$ is a
po-neighborhood, then any subset of $W_X(x)$ which is a neighborhood of $x$
is also a po-neighborhood\thinspace;
hence one can assume that $W(x)\subset U$ for some
$U\in\U$ and hence that the partial order on $W_X(x)$ is induced
from $U$. The po-neighborhoods satisfying this extra condition are
called {\em small po-neighborhoods} and they give a neighborhood basis
for the topology on $X$, since the intersection of two small
po-neighborhoods are again a small po-neighborhood. Moreover, the
covering by the small po-neighborhoods defines the local partial
order.

Let $X$ be a local po-space: it is defined by a covering
$({\cal  U},(\leq_U)_{U \in {\cal U}})$ of open sub-po-spaces of $X$ together
with $(W_X(x),\leq_{W_X(x)})$, for all $x \in X$, the local
neighborhood and the corresponding partial order. We now only consider
any of its small representatives in its equivalence class of
local po-structures (we still call $X=({\cal U}, (\leq_U)_{U \in
  {\cal U}})$).  $\vec{D}^n$ is a local po-space, which is actually a
(global) po-space. So its covering is $(\vec{D}^n,\leq_{\vec{D}^n})$
with corresponding $(W_{\vec{D}^n}(x)=\vec{D}^n,
\leq_{W_{\vec{D}^n}(x)}=\leq_{\vec{D}^n})$.

Let $f:\vec{S}^{n-1}\longrightarrow X$ be an attaching map\footnote{We consider
one attaching map at a time only, the argument easily transposes to any
number of attaching maps.} of a
globular  $n$-cell $\vec{e}^n$. We construct the topological
space $Z= \vec{D}^n \sqcup_{f} X $ as defined by Lemma~\ref{base}.
Let $\Phi: \vec{D}^n \sqcup X \rightarrow \vec{D}^n \sqcup_f X$ be
the canonical surjective map. We have a commutative diagram in the
category of topological spaces:
$$
\begin{diagram}
\vec{S}^{n-1} & \rTo^f & X \\
\dTo^{i} & & \dTo_{\Phi_1} \\
\vec{D}^{n} & \rTo_{\Phi_2} & Z \\
\end{diagram}
$$
where $i$ is the inclusion map and $\Phi_1$, $\Phi_2$ are defined
by the push-out diagram. Of course, $\Phi_1$ is injective since
$i$ is injective. We identify $\Phi_1$ with $\Phi$ restricted to
$X$ and also identify $\Phi_2$ with $\Phi$ since it is the
composition of the inclusion map from$\vec{D}^{n}$ to $\vec{D}^n \sqcup X$
with $\Phi$.

As $\vec{S}^{n-1}$ is compact, by Lemma~\ref{base}, point 3 and 4, we know
that $\Phi$ is a closed map and $Z$ is Hausdorff (this holds
true by induction again). Therefore $f(\vec{S}^{n-1})$ is closed since it is compact.
Thus by point 3 of Lemma~\ref{base}, $X \backslash \vec{S}^{n-1}$ can
be identified  with the open subset $\Phi(X \backslash f(\vec{S}^{n-1}))$ of
$Z$ and $f(\vec{S}^{n-1})$ with the closed subset
$\Phi(f(\vec{S}^{n-1}))$ of $Z$.

Similarly, $\vec{S}^{n-1}$ is a closed subset of $\vec{D}^n$ so by point 2 of
Lemma~\ref{base}, $\vec{D}^n \backslash f(\vec{S}^{n-1})$ can be identified
with the open subset $\Phi(\vec{D}^n \backslash f(\vec{S}^{n-1}))$ of $Z$ and
$X$ can be identified with the closed subset $\Phi(X)$ of $Z$. We use these
identifications below.

Take now $z \in Z$; we are going to construct a neighborhood
$U_z$ of $z$ in $Z$ together with a local po-structure on $U_z$,
making $Z$ into a local po-space with the local po-structure
$(U_z,\prec_{U_z})_{z\in Z}$\thinspace:
\begin{itemize}
\item [(1)] Suppose $z \in \vec{D}^n \backslash f(\vec{S}^{n-1})$
(see Figure~\ref{case1}). We define
$U_z=\vec{D}^n \backslash f(\vec{S}^{n-1})$ that we noticed is identified with
an open subset of $Z$,
and a binary
relation $\prec_{U_z}$
on $U_z$ such that $u \prec_{U_z} v$ if $u \leq_{\vec{D}^n} v$. $\prec_{U_z}$ is obviously
a partial order.
\begin{figure}
\begin{center}
\epsfig{file=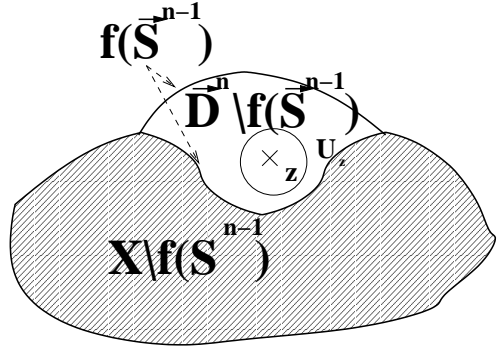,width=7cm,clip=}
\end{center}
\caption{First case for the definition of $(U_z,\prec_{U_z})$.}
\label{case1}
\end{figure}
\item [(2)] Suppose $z \in X \backslash f(\vec{S}^{n-1})$
(see Figure~\ref{case2}). We have noticed that
$X \backslash f(\vec{S}^{n-1})$ can be identified with an open subset of $Z$.
$W_X(z)$ is an open subset of $\vec{D}^n \sqcup X$ containing $\Phi^{-1}(\{z\})=\{z\}$
since $z$ is in $X$ and $\Phi$ is injective on this
part. Therefore, by Lemma~\ref{ferm}, there exists $U_z$ open of $Z$ containing
$\{z\}$ such that $\Phi^{-1}(U_z) \subseteq W_X(z)$.
We define the partial ordering $\prec_{U_z}$ to be the same as $\leq_{W_X(z)}$ on $U_z$.
\begin{figure}
\begin{center}
\includegraphics[width=7cm]{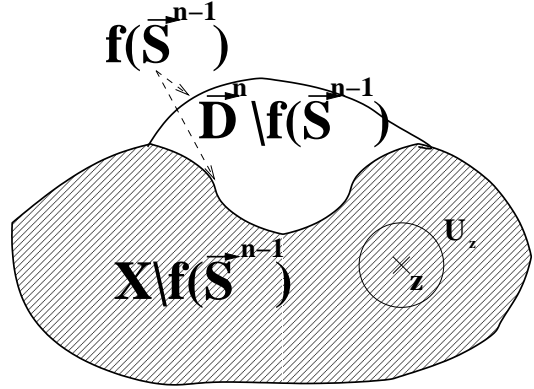}
\end{center}
\caption{Second case for the definition of $(U_z,\prec_{U_z})$.}
\label{case2}
\end{figure}
\item [(3)] The only remaining possibility is that $z \in f(\vec{S}^{n-1})$
(see Figure~\ref{case3}).
\begin{figure}
\begin{center}
\epsfig{file=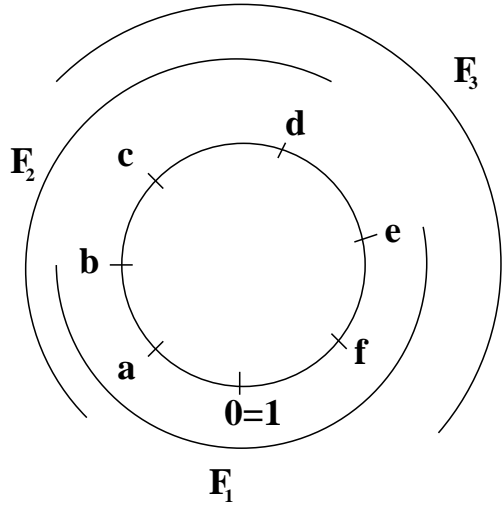,width=7cm,clip=}
\end{center}
\caption{The subdivision of an oriented circle.}
\label{subdDn}
\end{figure}
\begin{figure}
\begin{center}
\epsfig{file=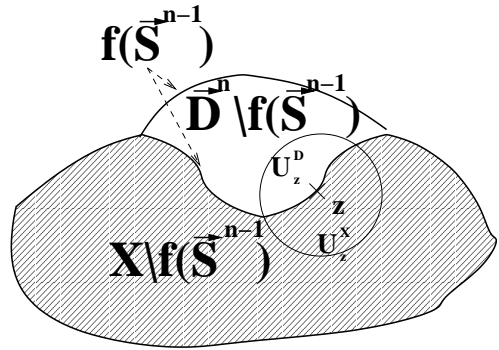,width=7cm,clip=}
\end{center}
\caption{Third case for the definition of $(U_z,\prec_{U_z})$.}
\label{case3}
\end{figure}
Let us first subdivide the segment $[0,1]$; take six elements of $]0,1[$
$0 < a < b < c < d < e < f < 1$. We let (see Figure~\ref{subdDn}),
\begin{itemize}
\item $F_1=[e,1] \cup [0,b]$, with partial order $\leq_{F_1}$ defined by,
for $x \in [e,1]$ and $y \in [e,1]$ or $x \in [0,b]$ and $y \in [0,b]$,
it is the usual partial order induced by $[0,1]$ and for $x \in [e,1]$
and $y \in [0,b]$, $x \leq_{F_1} y$.
\item $F_2=[a,d]$, with the usual partial order.
\item $F_3=[c,f]$, with the usual partial order.
\end{itemize}

We notice that if we identify $0$ with $1$, $(\intr{F_1},\leq_{F_1})$,
$(\intr{F_2},\leq_{F_2})$ and $(\intr{F_3},\leq_{F_3})$ is a small
local po-structure on the circle
 and the canonical surjection from the po-space
$\vec{I}$ to this local po-space is a morphism of local po-spaces.

We define $\vec{D}^n_i=\{(x,t) \mid x \in D^{n-1}, t \in F_i\}$
(similarly $\vec{S}^{n-1}_i=\{(x,t) \mid x \in S^{n-2}, t \in F_i\}$)
closed subset of $\vec{D}^n$.
The partial orders $\leq_{F_i}$ induce partial orders $\leq_{\vec{D}^n_i}$
on $\vec{D}^n_i$.

As $X$ is locally compact, we can
find $W_z$ a closed neighborhood of $z$ contained in $W_X(z)$.

Consider the composite map $\Phi_{i}$:

\[
\xymatrix{
\vec{D}^n_{i} \sqcup W_z \ar[dr]_{\Phi_i} \ar[r]^{\subseteq} & \vec{D}^n \sqcup X
\ar[d]_{\Phi} \\
& \vec{D}^n \sqcup_f X \\ }\]

It is a closed continuous map as a composition of two closed
continuous maps. There exists (a non-necessarily unique) $(w,t) \in
S^{n-2} \p \vI$ such that $f(w,t)=z$. Necessarily, $t$ belongs to some
$\intr{F_{i_z}}$.  We have $$\Phi_{i_z}^{-1}(\{z\}) \subseteq
\vec{D}^n_{i_z} \sqcup W_z$$
thus by Lemma~\ref{ferm} there exists
an open neighborhood $U_z$ of $z$ such that
$$ \Phi_{i_z}^{-1}(U_z) \subseteq
\vec{D}^n_{i_z} \sqcup W_z$$

Let $U_z^X$ be the
subset $U_z \cap \Phi(X)$ of $Z$ and $U_z^D$ be the subset
$U_z \cap (\Phi(\vec{D}^n \backslash f(\vec{S}^{n-1})))$ of $Z$.
This is a
partition of $U_z$.
Notice that we can identify elements of
$U_z^D$ with elements of $\vec{D}^n\backslash \vec{S}^{n-1}$
and elements of $U_z^X$ with
elements of $X$. By construction, $U_z^D \subseteq \vec{D}^n_{i_z}
\backslash \vec{S}^{n-1}_{i_z}$.
We now define a binary relation $\prec_{U_z}$ on
$U_z$ as follows:
\begin{itemize}
\item for $u, v \in U^X_z$, $u \prec_{U_z} v$ if $u \leq_{W_X(z)} v$,
\item for $u, v \in U^D_z$,
$u \prec_{U_z} v$ if $u \leq_{\vec{D}^n_{i_z}}v$,
\item for $u \in U^X_z$ and $v \in U^D_z$,
\begin{itemize}
\item if $i_z=1$, $u \prec_{U_z} v$ if $u \leq_{W_X(z)} f(\ei)$ and
$0 < t(v) < a$, ($t(v)$ is the unique parameter in $F_1$ such
that $v=(w,t(v))$ for some $w$),
\item if $i_z=2$, $u \prec_{U_z} v$  if $u \leq_{W_X(z)} f(\ei)$,
\item if $i_z=3$ we can never have $u \prec_{U_z} v$.
\end{itemize}
\item for $u \in U^D_z$ and $v \in U^X_z$,
\begin{itemize}
\item if $i_z=1$, $u \prec_{U_z} v$ if $f(\ef) \leq_{W_X(z)} v$ and
$b < t(u) < 1$,
\item if $i_z=2$ we can never have $u \prec_{U_z} v$,
\item if $i_z=3$, $u \prec_{U_z} v$ if $f(\ef) \leq_{W_X(z)} v$.
\end{itemize}
\end{itemize}

This defines a partial order indeed. Reflexivity and transitivity are obvious.
We
now check antisymmetry. Let $u$ and $v$ be two elements of $U_z$ such
that $u \prec_{U_z} v$ and $v \prec_{U_z} u$. If $u$ and $v$ both
belong to $U^X_z$ or $U^D_z$ it is obvious that this implies $u=v$,
since the relation $\prec_{U_z}$ coincide with one of the partial
orders $\leq_{W_X(z)}$ or $\leq_{\vec{D}^n_i}$ in that case. Suppose
$u \in U^X_z$, $v \in U^D_z$ with,
\begin{itemize}
\item $i_z=1$, we have by definition $u \leq_{W_X(z)} f(\ei)$ and $0 <
t(v) < b$ and $f(\ef) \leq_{W_X(z)} u$ and $e < t(v) < 1$, which is
of course impossible,
\item $i_z=2$, $v \prec_{U_z} u$ is impossible by definition,
\item $i_z=3$, $u \prec_{U_z} v$ is impossible by definition.
\end{itemize}

\end{itemize}

It follows that $(U_z,\prec_{U_z})_{z\in Z}$ defines a small
local po-structure since by construction, for $z\neq z'$, the
partial orders $\prec_{U_z}$ and $\prec_{U_{z'}}$ coincide on
the intersection $U_z\cap U_{z'}$ (if non-empty). It then suffices to set $W_Z(z):=U_z$. \epf

\bth The previous embedding induces a functor from the category
of globular CW-complexes to that of local po-spaces. \eth

\bpf Let $X$ and $Y$ be two globular CW-complexes and $f: X \rightarrow Y$ be a
morphism of CW-complexes. The globular cellular decomposition of $X$ yields
a set of characteristic maps $\phi_\alpha:\vec{D}^{n_\alpha}\rightarrow X$
satisfying\thinspace:
\begin{enumerate}
\item The mapping
$\phi_\alpha\restriction_{\vec{D}^{n_\alpha}-\vec{S}^{n_\alpha-1}}$
induces an homeomorphism from $\vec{e}^{n_\alpha}$ to its image.
\item All the previous globular cells are disjoint and their union
gives back $X$.
\item A subset of $X$ is closed if and only if it meets the closure
of each globular cells of $X$ in a closed set.
\end{enumerate}
where $\alpha$ runs over a well-ordered set of indexes $\kappa$ (one can suppose
that $\kappa$ is a finite or transfinite cardinal). One can suppose that the mapping
$\alpha \mapsto n_\alpha$ is non-decreasing. Let $X^{[-1]}=\emptyset$.
Let $\beta$ be an ordinal
with $\beta \leq \kappa$. If $\beta$ is a limit ordinal, let
$X^{[\beta]}=\varinjlim_{\alpha<\beta} X^{[\alpha]}$. If $\beta=\gamma+1$ for some
ordinal $\gamma$, then let $X^{[\beta]}= \vec{D}^{n_\beta} \sqcup_{\phi_\beta} X^{[\gamma]}$.
Notice that $X^{[\gamma]}$ is closed in $X^{[\beta]}$.

We are going to prove by transfinite induction on $\beta$ the
statement $P(\beta)$\thinspace: \textit{For any globular CW-complex $X$ and for any set of
characteristic maps $\phi_\alpha:\vec{D}^{n_\alpha}\rightarrow X$ as above,
a morphism of globular CW-complexes from
$X$ to $Y$ induces a morphism of local po-spaces from $X^{[\beta]}$ to
$Y$.}

Necessarily the equality $n_0=0$ holds therefore $P(0)$ is true. Now
let us suppose that $P(\alpha)$ holds for $\alpha<\beta$ and some
$\beta \geq 1$. We want to check that $P(\beta)$ then holds as well.
If $\beta=1$, then $X^{[\beta]}$ is either the two-point discrete
space, or a loop. So $P(1)$ holds. So let us suppose $\beta \geq 2$.

There are two mutually exclusive cases\thinspace:
\begin{enumerate}
\item The case where $\beta$ is a limit ordinal. Let $x\in X^{[\beta]}$.
Then $x\in X^{[\alpha]}$ for some $\alpha<\beta$ and the induction
hypothesis can be applied\thinspace; the result follows from the fact that the direct
limit is endowed with the weak topology.
\item The case where $\beta=\gamma+1$ for some cardinal $\gamma$. Then
$X^{[\beta]}= \vec{D}^{n_\beta} \sqcup_{\phi_\beta} X^{[\gamma]}$ with
the above notations. With the notation and identification as in
the proof of Theorem~\ref{CWpo}, one has
three mutually exclusive cases\thinspace:
\begin{itemize}
\item $x\in X^{[\gamma]}\backslash \phi_\beta(\vec{S}^{n_\beta-1})$\thinspace: in this case,
the induction hypothesis can be applied\thinspace;
\item $x\in \vec{D}^{n_\beta}\backslash \phi_\beta(\vec{S}^{n_\beta-1})$\thinspace: let $W_Y(f(x))$
be a po-neighborhood of $f(x)$ in $Y$\thinspace; then $f^{-1}(W_Y(f(x)))$ is an open of
$\vec{D}^{n_\beta}$\thinspace; there exists a basis of $\vec{D}^{n_\beta}$ by global po-spaces so
there exists a po-neighborhood $W_x$ of $x$ in $\vec{D}^{n_\beta}$ such that
$W_x\subset f^{-1}(W_Y(f(x)))$\thinspace;
\item $x\in \phi_\beta(\vec{S}^{n_\beta-1})$\thinspace: let $\Phi_\beta$ be the canonical closed map
from $\vec{D}^{n_\beta} \sqcup X^{[\gamma]}$ to $X^{[\beta]}$\thinspace; by induction hypothesis,
$f\circ \Phi_\beta\restriction_{X^{[\gamma]}}:X^{[\gamma]}\rightarrow Y$ is a morphism
of po-spaces\thinspace; therefore there exists a po-neighborhood $W_{f(x)}$ of $f(x)$ in $Y$ and
a po-neighborhood $W_{x}$ of $x$ in $X^{[\gamma]}$ such that
\[ W_x \subset \left(f\circ \Phi_\beta\restriction_{X^{[\gamma]}}\right)^{-1}\left(W_{f(x)}\right)\]
So $(\Phi_\beta)^{-1}(\{x\})\in \vec{D}^{n_\beta} \sqcup W_x$ and by Lemma~\ref{ferm},
there exists an open $V_x$ of $X^{[\gamma]}$ such that
$(\Phi_\beta)^{-1}(V_x)\in \vec{D}^{n_\beta} \sqcup W_x$. Then let us considering the $U_x$
of the proof of Theorem~\ref{CWpo}. Since $f$ is continuous, $f^{-1}(W_{f(x)})$ is open and
\[\emptyset \neq V_x \cap U_x \cap f^{-1}(W_{f(x)}) \subset f^{-1}(W_{f(x)})\]
\end{itemize}

\end{enumerate}

\epf

We now prove an interesting tool for constructing globular complexes.

\bth Let $Z$ be a compact local po-space, let $Y$ be a closed subset  of
$Z$, and let $\e$ be a globular $n$-cell in $Z$ with $\e\cap Y=\emptyset$.
Suppose there
exists a relative isomorphism\footnote{Meaning that $\Phi$
is an isomorphism of globular CW-complexes from $\vec{D}^n$ to $\e \cup Y$
such that it restricts to an isomorphism of globular CW-complexes from
$\vec{S}^{n-1} \subseteq \vec{D}^n$ to $Y \subseteq \e \cup Y$.}
of globular CW-complexes
$\xymatrix@1{\Phi:(\vec{D}^n,\vec{S}^{n-1})\fr{} & (\e \cup Y,Y)}$.
Set $f=\Phi|\vec{S}^{n-1}$.
then the obvious map (induced by $\Phi$ and by $Id_{Y}$)
\[\Psi:\xymatrix@1{{Y_f=\vec{D}^n \sqcup_f Y}\fr{} & {\e\cup Y}}\]
is an isomorphism of local po-spaces.
\eth

\bpf
The map $\Psi$ is clearly bijective.
Let $g$ be the canonical map from $\vec{D}^n \sqcup Y$ to
$\vec{D}^n \sqcup_f Y$ and let
$\Omega$ be an open subset of $\e\cup Y$. Then
$g^{-1}\Psi^{-1}(\Omega)=\Phi^{-1}(\Omega\cap \e)\sqcup (\Omega\cap Y)$
is an open subset of $\vec{D}^n \sqcup Y$, therefore $\Psi$ is
continuous.  So $\e \cup Y$ is compact and therefore $\Psi$ is an
homeomorphism. The fact that $\Psi$ preserves the structure of local
po-spaces is obvious.
\epf

\subsection{Globular CW-complex and precubical set}

\label{appA}
\label{A}

We are going to show that in fact, there is a geometric realization
functor (which should be homotopically equivalent to the former one composed
with the realization of precubical sets as local po-spaces of
\cite{LFEGMRAlgebraic} in
some sense) transforming a precubical set into a globular
CW-complex. We first need a few (classical) remarks.

\bd\label{def_cubique}\cite{Brown_cube} \cite{cube} A
\textit{precubical set} (or HDA) consists of a family of sets
$(M_n)_{n\geqslant 0}$ and of a family of face maps
$\xymatrix@1{{M_n}\fr{\de_i^\alpha} &{M_{n-1}}}$ for
$\alpha\in\{0,1\}$ and $1\leq i \leq n$
which satisfies the following axiom (called
sometime the cube axiom)\thinspace:
\begin{center}
$\de_i^\alpha \de_j^\beta = \de_{j-1}^\beta \de_i^\alpha$ for all
$1\leq i<j\leqslant n$ and $\alpha,\beta\in\{0,1\}$.
\end{center}
\ed

If $M$ is a precubical set, the elements of $M_n$ are called the
$n$-cubes. An element of $M_n$ is of dimension $n$. The elements
of $M_0$ (resp. $M_1$) can be called the \textit{vertices} (resp.
the \textit{arrows}) of $K$.

Let $M$ and $N$ be two HDA, and $f$ a family $f_{n}: M_n
\rightarrow N_n$ of functions. The family $f$ is a morphism of
HDA if and only if $ f_{n} \circ \de^\alpha_i = \de^\alpha_i
\circ f_{n+1} $ for all $n \in \N$ and $0 \leq i \leq n$.  HDA
together with these morphisms form a category which we denote by
$\square {\bf Set}$. Conventionally, this category can be identified
with the set-valued pre-sheaves of some small finite free category
$\square$, and therefore it is cocomplete.

The small category $\square$ can be described as a category whose
objects are $[n]$, where $n \in \N$ and whose morphisms are
generated by, $$
\CellSize=1.75em\begin{diagram}
[n-1] & \pile{\rTo^{\delta^0_i} \\ \rTo_{\delta^1_j}} & [n] \\
\end{diagram}
$$
with $1\leq i,j\leq n$ and satisfying the opposite of the cube axiom,
i.e. $ \delta_j^\beta\delta_i^\alpha =  \delta_i^\alpha \delta_{j-1}^\beta$ for all
$1\leq i<j\leqslant n$ and $\alpha,\beta\in\{0,1\}$.

There is a truncation functor $T_n: \square {\bf Set} \rightarrow \square {\bf Set}$ defined by,
$T_n(M)_{m}=M_{m}$ if $m \leq n$ and $T_n(M)_{m}=\emptyset$ if $m > n$.

Now, let $D_{[n]}$ be the representable functor from $\square$ to $\bf Set$
with $D_{[n]}([p])=\square {\bf Set}([p],[n])$. We define
the singular $n$-cubes of a HDA $M$ to be any morphism $\sigma: D_{[n]} \rightarrow M$.

\begin{lem}
The set of singular $n$-cubes of a HDA $M$ is in one-to-one
correspondence with $M_n$. The unique singular $n$-cube
corresponding to a $n$-cube $x \in M_n$ is denoted by $\sigma_x:
D_{[n]} \rightarrow M$. It is the unique singular $n$-cube
$\sigma$ such that $\sigma(Id_{[n]})=x$.
\end{lem}

\bpf
The proof goes by
Yoneda's lemma. 
\epf

There is only one morphism in $\square$ from a given $[n]$ to itself, the
identity of $[n]$,
hence $D_{[n]} \backslash \{Id\}$ is a functor which has only as non-empty
values the $D_{[n]}([p])$ with $p < n$ (``it is of dimension $n-1$'').
We write $\de{D}_{[n]}$ for this functor. For $\sigma$ a natural transformation
from $D_{[n]}$ to $M$, we write $\de{\sigma}$ for its restriction to $\de{D}_{[n]}$.

\bp
Let $M$ be a HDA. The following diagram is co-cartesian
(for $n \in \N$),
$$ \CellSize=1.75em \begin{diagram}
\coprod_{x \in M_{n+1}} \de{D}_{[n+1]} &
\rTo^{\bigsqcup_{x \in M_{n+1}} \de{\sigma}_x} & T_n(M) \\
\dTo^{\subseteq} & & \dTo^{\subseteq} \\
\coprod_{x \in M_{n+1}} D_{[n+1]} &
\rTo^{\bigsqcup_{x \in M_{n+1}} \sigma_x}
& T_{n+1}(M) \\
\end{diagram} $$
where $\de{D}_{[n+1]}=T_n(D_{[n+1]})$ and $\de{\sigma}_x= {\sigma_x}_{|\de{D}_{[n+1]}}$.
\ep

\bpf We mimic the proof of \cite{PGMZCalculus}: it suffices to
prove that the diagram below (in the category of sets) is
cocartesian for all $p \leq n+1$,
$$ \CellSize=1.75em \begin{diagram}
\coprod_{x \in M_{n+1}} (\de{D}_{[n+1]})_p &
\rTo^{\bigsqcup_{x \in M_{n+1}} (\de{\sigma}_x)_p} & (T_n(M))_p \\
\dTo^{\subseteq} & & \dTo^{\subseteq} \\
\coprod_{x \in M_{n+1}} (D_{[n+1]})_p &
\rTo^{\bigsqcup_{x \in M_{n+1}} (\sigma_x)_p}
& (T_{n+1}(M))_p \\
\end{diagram} $$
since colimits (hence push-outs) are taken point-wise in a functor
category into $\bf Set$.

For all $p < n+1$, the inclusions are in fact bijections, and the
diagram is then obviously cocartesian.

For $p=n+1$, the complement of $\bigsqcup_{x \in M_{n+1}} (\de{D}_{[n+1]})_p$
in $\bigsqcup_{x \in M_{n+1}} (D_{[n+1]})_p$ is the set of copies of cubes
$Id_{[n+1]}$, one for each cube of $M_{n+1}$. This means that
the map $\bigsqcup_{x \in M_{n+1}} (\sigma_x)_p$ induces a bijection
from the complement of $\bigsqcup_{x \in M_{n+1}} (\de{D}_{[n+1]})_p$
onto the complement of $(T_n(M))_p$. This implies that the
diagram is cocartesian for $p=n+1$ as well.
\epf

This lemma states that the truncation of dimension $n+1$ of a
HDA $M$ is obtained from the truncation of dimension
$n$ of $M$ by attaching some standard $(n+1)$-cubes $D_{[n+1]}$ along the
boundary $\de{D}_{[n+1]}$ of $(n+1)$-dimensional holes.
In fact, any precubical set $M$ is the direct limit of the diagram
consisting of all inclusions $T_{n-1}(M) \hookrightarrow T_n(M)$,
hence is also the direct limit of the diagram consisting of all the cocartesian
squares above.

We are quite close
to the globular CW-complex definition. What we need now is the (classical)
notion of geometric realization.
Let $\square_n$ be the standard cube in $\R^n$ ($n \geq 0$), $$
\square_n = \{(t_1,\ldots,t_n) | \forall i, 0 \leq t_i \leq 1 \} $$ $$
\square_{0} = \{0\} $$ and let $\delta^k_i:\square_{n-1}\to\square_n$, $1 \leq i \leq n$, $k=1,2$, be the
continuous functions ($n \geq 1$), $$
\begin{diagram}
\square_n & \lTo^{\delta^0_i} & \square_{n-1} \\
\uTo^{\delta^1_i} & & \\
\square_{n-1} & & \\
\end{diagram}
$$ defined by, $$
\begin{array}{rcl}
\delta^k_i (t_1,\ldots,t_{n-1}) & = & (t_1,\ldots,t_{i-1},k,t_i,\ldots,t_{n-1})
\\
\end{array}
$$

Consider now, for a precubical set $M$, the set ${\bf
R}(M)=\bigsqcup_{n} M_n\times\square_{n}$.  The sets $M_n$ have the
discrete topology and $\square_{n}$ is endowed with the topology
induced as a subset of ${\R}^n$  with the standard topology thus
${\bf R}(M)$ is a topological space with the disjoint sum topology.\\
Let $\equiv$ be the equivalence relation induced by the
identities:

$$\forall k,i,n, \forall x \in M_{n+1}, \forall t \in \square_{n}, n
\geq 0, (\partial^k_i(x),t) \equiv (x,\delta^k_i(t))$$

Let $|M| = {\bf R}(M)/\equiv$ have the quotient topology. The topological space $| M|$ is called the {\em geometric realization} of $M$. This actually
yields  a functor from $\square{\bf Set}$ to $\top$ by setting
for $f: X \rightarrow Y$ a morphism in $\square{\bf Set}$,
${\bf R}(f): {\bf R}(X) \longrightarrow {\bf R}(Y)$ by:
${\bf R}(f)((x,t))=(f(x),t)$. This functor commutes with colimits since it is a left
adjoint functor (the right adjoint being a singular cube functor, see
\cite{LFEGMRAlgebraic}).

Thus, the geometric realization of a precubical set $M$ is the direct
limit of the diagram:
$$ \CellSize=1.75em \begin{diagram}
\coprod_{x \in M_{n+1}} \de{\square}_{n+1} &
\rTo^{\bigsqcup_{x \in M_{n+1}} | \de{\sigma}_x|} & | T_n(M)| \\
\dTo^{\subseteq} & & \dTo^{\subseteq} \\
\coprod_{x \in M_{n+1}} \square_{n+1} &
\rTo^{\bigsqcup_{x \in M_{n+1}} | \sigma_x |}
& | T_{n+1}(M)| \\
\end{diagram} $$
since it is easily shown that,
\begin{itemize}
\item $| D_{[n+1]} |$ is homeomorphic to $\square_{n+1}$,
\item the inclusion of $\de{D}_{[n+1]}$
into $D_{[n+1]}$ induces an homeomorphism between $| \de{D}_{[n+1]} |$
onto the boundary $\de{\square}_{n+1}$ of the standard $(n+1)$-cube.
\end{itemize}
Obviously, each $\sigma_x$ induces a homeomorphism from $\de{\square}_{n+1}$
onto a connected component of $| T_{n+1}(M)| \backslash | T_n(M) |$,
which is homeomorphic to the interior of an $(n+1)$-cube, and to $e^{n+1}$.
This shows that at least, $| M |$ is a CW-complex. We are now going to
show that this direct limit can be slightly transformed so as to produce
a globular CW-complex.

Let $C_t=\{ (t_1,\dots ,t_n)\in \square _n|t_1+\dots +t_n=t\} ,\; 0\le
t\le n$.  As the intersection of $\square _n$ with a hyperplane, $C_t$ is
a convex polyhedron. It is spanned by (the convex hull of) extremal
points all of whose coordinates but one are either $0$ or $1$. More
precisely, extremal points have $\lfloor t\rfloor $ 1-coordinates,
$(n-\lfloor t\rfloor -1)$ 0-coordinates and a single coordinate is
given by $t-\lfloor t\rfloor $.  ($\lfloor \, \rfloor $ denotes Gauss
brackets, i.e., $\lfloor t\rfloor $ is the largest integer $\le
t$.)\footnote{This generalizes the picture in dimension 3, where $C_t$
  is either a triangle or a hexagon (for $1<t<2$).}
In particular, the point $(\frac{t}{n},\dots ,\frac{t}{n})$ is an
interior point in $C_t$ for $0<t<n$. For $t=0$ and for $t=n$, the
polyhedron $C_t$ degenerates to a single point.

Let $h:\mb{R}^n\to \mb{R}^{n-1}$ denote the linear map given
by
$$h(t_1,\dots ,t_n)=(t_2-\frac{t_1+\dots +t_n}{n},\dots
,t_n-\frac{t_1+\dots +t_n}{n}).$$
Its kernel $\ker (h)$ is equal to
the diagonal $\Delta \subset \mb{R}^n$.  Let $D_t=h(C_t)\subset
\mb{R}^{n-1},\; 0\le t \le n$. As the image of a convex polyhedron
under an injective affine map (on the hyperplane $\sum t_i=t$), $D_t$
is itself a convex polyhedron containing the origin
$\mb{0}=h(\frac{t}{n},\dots ,\frac{t}{n})$ in its interior for $0<t<n$;
$D_t$ degenerates to a single point exactly for $t=0$ and $t=n$.

Define $m:\square _n\setminus \Delta \to \mb{R}$ by $m(t_1,\dots
,t_n)=\min _i\max \{ m_i\ge 0|\; 0\le
\frac{t}{n}+m_i(t_i-\frac{t}{n})\le 1\} $. ($m$ is continuous, but not
well-defined on the diagonal $\Delta \subset \square _n$; on the ray from
the diagonal point in $C_t$ through $\mb{t}=(t_1,\dots ,t_n)$ to the
boundary of $C_t$, the number $\frac{1}{m(\mb{t})}$ measures
the ratio of the segment between the diagonal and $\mb{t}$). A
homeomorphism $H:\square _n\to Glob(D^{n-1})$ can be defined as the
composition of the quotient map $q:I\times D^{n-1}\to Glob(D^{n-1})$
with the following map $\bar{H}:\square _n\to I\times D^{n-1}$:
$$\bar{H}(t_1,\dots ,t_n)=\left\{ 
  \begin{array}{ll}
(\frac{t_1+\dots +t_n}{n}, \frac{h(t_1,\dots
  ,t_n)}{m(t_1,\dots ,t_n)||h(t_1,\dots ,t_n)||}) &  (t_1,\dots
,t_n)\not\in \Delta ,\\
(\frac{t_1+\dots
  +t_n}{n},0,\dots ,0) & (t_1,\dots ,t_n)\in \Delta .
  \end{array}\right. $$
The function $m\equiv 1$ on $\partial \square _n\setminus \{
\mb{0},\mb{1}\} $, whence $H(\partial \square _n\setminus \{
\mb{0},\mb{1}\})=int(I)\times S^{n-2}$ and $H(\square _n\setminus \{
\mb{0},\mb{1}\})=int(I)\times D^{n-1}$. The restriction of $\bar{H}$
to $\square _n\setminus \{ \mb{0},\mb{1}\} $ is a continuous bijection
onto its image. Remark, that $H$ -- but not $\bar{H}$ -- is
continuous in the corners points, $(0,\dots ,0)$ and $(1,\dots 1)$. As
a continuous bijection from the compact space $\square _n$ into the
Hausdorff space $Glob(D^{n-1})$, the map $H$ is a homeomorphism.

Define a partial order on $\square_n$
by $x \leq_{gl} y$ if and only if $H(x) \leq H(y)$ (using the partial
order on $Glob(D^{n-1})$). Then by definition, $(\square_n,\leq_{gl})$ is isomorphic
as a po-space to $\vec{D}^n$ (through $H$). Notice also that
$H^{-1}(\de{\vec{D}}^{n})=\de{\square}_n$. We are now ready for the
construction of a globular CW-complex out of $M$.

\begin{itemize}
\item We start with $X^0=\coprod_{x \in M_0} \square_0$.
\item We form inductively the $n$-skeleton $X^n$ from $X^{n-1}$, that
we prove (again by induction on $n$) to be homeomorphic to $| T_{n-1}(M)|$,
by attaching globular $n$-cells $\vec{e}^n_{\sigma_x}$ via maps
$\phi_{\sigma_x}:\vec{S}^{n-1}\longrightarrow X^{n-1}$, where $\sigma_x$ is any
singular $n$-cube of $M$. The attaching map is defined as the composite:
$$\begin{diagram}
\phi_{\sigma_x}: & \vec{S}^{n-1} & \rTo^{H^{-1}} & \de{\square}_n &
\rTo^{| \de{\sigma}_x |} & | T_{n-1}(M) | \sim X^{n-1}
\end{diagram}
$$
What will remain to be shown
is that this is non-decreasing and
non-contracting.
\item
  $X=\bigcup_n X^n$ with the weak topology (the direct limit of the diagram
composed of the attaching maps).
\end{itemize}

We now check that the attaching maps are non-decreasing and non-contracting.
First of all, for any $x \in M_n$, we consider $\sigma_x: D_{[n]} \rightarrow  M$,
the unique morphism of precubical sets such that
$\sigma(Id_{D_{[n]}})=x \in M_n$.
We have to check that
if $y \leq_{gl} z$ in $(Id_{D_{[n]}},\de{\square}_n)$, then
$| \sigma_x |(x) \leq_{gl} | \sigma_x |(z)$
in $(x,\square_n) \in M$.
All points of $| D_{[n]} |$ have a representative in
$(Id_{D_{[n]}},\square_n)$, i.e. can be written as
$(Id_{D_{[n]}},t_1,\cdots,t_n)$ with $0 \leq t_i \leq 1$ (for all $i$,
$1 \leq i \leq n$). Now,
$| \sigma_x |(Id_{D_{[n]}},t_1,\cdots,t_n)= (x,t_1,\cdots,t_n)$,
hence $| \sigma_x |$ preserves trivially the
partial order $\leq_{gl}$ of $\square_n$,
hence $| \sigma_x | \circ H^{-1}$ preserves it as well. 
Since
$H^{-1}(\de{\vec{D}}^n)=\de{\square}_n$,
$\de{\sigma}_x$ also preserves $\leq_{gl}$.

Now, take an execution path $\phi$ starting from $\ei$ (or arriving at
$\ef$) in $\vec{S}^{n-1}$, and suppose that $\phi_{\sigma_x}\circ
\phi$ is a constant path of $X^{n-1}$. Then $\sigma_x \circ
H^{-1} \circ \phi$ has constant coordinates in $(x,\square_n) \in
| T_n(M) |$, which means, since $\sigma_x$ acts as the identity
on these coordinates, that $\phi$ is a constant path of
$\vec{S}^{n-1}$.  Furthermore,
$\phi_{\sigma_x}(\ei)=\sigma_x(Id_{[n]},(0,\cdots,0))$ which is also
$\sigma_x(\delta^0_0 d^0_1 \cdots \delta^0_{n-1},\square_0)$, so is equal
to $(\delta^0_0 \delta^0_1 \cdots \delta^0_{n-1}(x),\square_0)$ which
belongs to $T_n(M)_0=M_0$.  Similarly, $\phi_\alpha(\ef)$ belongs to
$X^0$.

\bp The above construction induces a functor from the category of HDA
$\square{\bf Set}$ to the category $\diCW$ of globular CW-complexes. \ep

\bpf By definition, a morphism a precubical set sends a
$n$-cube to another $n$-cube. So the realization as globular
CW-complexes induce clearly a morphism of $\diCW$.  \epf

\begin{ex}\label{runningexample2} Consider the $2$-dimensional 
precubical set 
$sw=(sw_0,sw_1,sw_2)$ with (front face) boundary maps $\partial^0_j:
sw_{i+1} \rightarrow sw_{i}$ ($i=0, 1$, $0 \leq j \leq i$) and (back face) boundary
maps $\partial^1_j: sw_{i+1} \rightarrow sw_{i}$ ($i=0, 1$, $0 \leq j \leq i$)
defined as follows:
\begin{itemize}
\item 
The 0-cubes are $sw_0=\{0,1,2,3,4,5\}\p \{0,1,2,3,4,5\}$.
\item The 1-cubes are 
\beas
sw_1&=&\left\{((i,j),(i+1,j))\hbox{ for } (i,j)\in\{0,\dots,4\}\p \{0,\dots,5\}\right\}\\
&\cup& \left\{((i,j),(i,j+1))\hbox{ for } (i,j)\in\{0,\dots,5\}\p \{0,\dots,4\}\right\}\\
&\backslash & \left(
\{((2,2),(2,3)),((2,2),(3,2)), ((2,3),(3,3)),((3,2),(3,3))\}
\right)
\eeas
with boundary maps: $\partial^0_0((i,j),(i+1,j))=(i,j)$, $\partial^1_0((i,j),(i+1,j))=
(i+1,j)$.
\item The 2-cubes are 
\beas sw_2 & = & \{((i,j),(i+1,j+1)) \mid 0 \leq i, j \leq 4\} \\
& \backslash & \{(i,j)\in\{(2,1),(1,2),(2,2),(3,2),(2,3)\} \\
\eeas with boundary maps:
$\partial^0_0((i,j),(i+1,j+1))=((i,j),(i+1,j))$, 
$\partial^0_1((i,j),(i+1,j+1))=((i,j),(i,j+1))$,
$\partial^1_0((i,j),(i+1,j+1))=((i,j+1),(i+1,j+1))$,
$\partial^1_1((i,j),(i+1,j+1))=((i+1,j),(i+1,j+1))$.
\end{itemize}
This precubical set gives the semantics of the term $Pa.Pb.Va.VB \mid
Pb.Pa.Vb.Va$ as shown in \cite{LFEGMRAlgebraic}. We see that the realization
of this precubical set as defined above is exactly the globular CW-complex
of Example~\ref{runningexample}.
\end{ex}

\subsection{Globular CW-complex and d-space}

Marco Grandis defines in \cite{grandis1,grandis2} a notion of d-space to study
the geometry of directed spaces up to a form of dihomotopy. There is
a functor which associates to each globular CW-complex $X$ a d-space
$d_X$ as follows:
\begin{itemize}
\item the underlying topological space of $d_X$ is the underlying topological
space of $X$,
\item the set of paths is the union of the space of execution paths
of $X$ between any two points of the 0-skeleton with the set of all
points (seen as constant paths). This set of paths is topologized with
the Kelleyfication of the compact-open topology (see Section
\ref{orthogonal}). The composition operation is the obvious one. 
\end{itemize}

\section{Dihomotopy equivalence}
\label{ST}

As pointed out in \cite{ConcuToAlgTopo}, there are two types of
deformations of HDA which leave unchanged its computer-scientific
properties\thinspace: the spatial ones and the temporal ones. The aim of this
section is to define in a precise manner these notions. In other
terms, we are going to construct three categories whose isomorphism
classes of objects are respectively the globular CW-complexes modulo
spatial deformations, modulo temporal deformations and modulo both
kinds of deformations together.

We meet in this section a few set-theoretic problems which must be
treated carefully. So two universes $\U$ and $\V$ with $\U\in \V$ are
fixed. The sets are the elements of $\U$. One wants to construct
categories whose objects are sets and whose the collection of
morphisms between any pair of objects is a set as well. So by
\textit{category}, one must understand a $\V$-small category with
$\U$-small objects and $\U$-small homsets. By \textit{large category},
one must understand a category with $\V$-small objects, and $\V$-small
homsets whose set of objects is not $\V$-small \cite{MR96g:18001a}.

\subsection{S-homotopy equivalence}

Intuitively, the spatial deformations correspond to usual deformations
orthogonally to the direction of time. This is precisely what a
S-homotopy does.

\bd\label{s-dihomotopic}
Let $f$ and $g$ be two morphisms from the globular CW-complex $X$ to
the globular CW-complex $Y$. Then $f$ and $g$ are \textit{S-homotopic}
if there exists a continuous map $H$ from $X\p \I$ to $Y$
such that (writing $H_t=H(-,t)$),
\begin{itemize}
\item $H_0=f$, $H_1=g$ and for any
$t\in [0,1]$,
\item $H_t$ is a morphism of globular CW-complexes from $X$ to $Y$.
\end{itemize}
When this holds, we write $f\sim_{S} g$. The map $H$ is called a
S-homotopy from $f$ to $g$. This defines an equivalence relation on
the set of morphisms of globular CW-complexes from $X$ to $Y$. The
quotient set is denoted by $[X,Y]_{S}$.  \ed

For comparison purposes, the set of continuous maps up to homotopy
from $X$ to $Y$ will be denoted by $[X,Y]$ and the corresponding
equivalence relation by $\sim$.

If $X$ and $Y$ are two globular CW-complexes, a S-homotopy
$H:X\p \I\rightarrow Y$ without further precisions means that for every
$t\in \I$, $H_t=H(-,t)$ is a morphism of globular CW-complexes from $X$
to $Y$ and therefore that $H$ is a S-homotopy between $H_0$ and
$H_1$.

In particular, this means that given $x\in X^0$, the image of the map
$t\mapsto H(x,t)$ is included in the discrete set $X^0$, and therefore
that it is the constant map. Therefore, if $f$ and $g$ are two
S-homotopic morphisms of globular CW-complexes, then {\em $f$ and
$g$ will coincide on the $0$-skeleton of $X$}.

The latter remark leads us to introducing \textit{the cylinder functor} $I^S$
associated to the notion of S-homotopy. If $X$ is a CW-complex, let
$I^S Glob(X)\thinspace:= Glob(X\p {\I})$.  Now for any globular
CW-complex $X$, let us define $I^S X$ inductively on the globular
cellular decomposition of $X$ in the following manner\thinspace:
\begin{enumerate}
\item [1)] Let
$I^S(X)^0:=X^0$\thinspace;
\item [2)] Let us suppose the $n$-skeleton $I^S(X)^n$ defined
for $n\geq 0$\thinspace; For every $(n+1)$-dimensional globular cell
$(Glob(D^n),\phi:Glob(S^{n-1})\rightarrow X^n)$ of $X$, the globular
CW-complex $Glob(\linebreak[0]D^n\p \I)$ is attached to $I^S(X)^n$ by the attaching
map $I^S\phi:Glob(S^{n-1}\p \I)\rightarrow I^S(X)^n$.
\item [3)] Then the direct
limit $I^S X$ is a globular CW-complex.
\end{enumerate}

Topologically, $I^S X$ is the quotient of $X\p I$ by the
relations $\{x_0\}\p \I = \{\overline{x_0}\}$ for every $x_0\in X^0$
(since $X^0$ is discrete, this relation is closed and the quotient is
then still Hausdorff). Let $\epsilon_h$ be the morphism from $X$ to
$I^S(X)$ with $\epsilon_h(x)=(x,h)$ and $\sigma$ be the canonical map
from $X\p \I$ to $I^S(X)$. Then,

\bp Let $f$ and $g$ be two morphisms of globular CW-complexes from $X$
to $Y$. Then $f$ and $g$ are S-homotopic if and only if
there exists a morphism of globular CW-complex $\overline{H}$ from $I^S(X)$ to
$Y$ such that $\overline{H}\circ \epsilon_0=f$ and $\overline{H}\circ \epsilon_1=g$.
\ep

\bpf If such a
$\overline{H}$ exists, then $\overline{H}\circ \sigma$ is a S-homotopy
from $f$ to $g$. Reciprocally, if $H$ is a S-homotopy
from $f$ to $g$, then the map $t\mapsto H(x_0,t)$ is
constant for any $x_0\in X^0$. Therefore $H$ factorizes by $\sigma$, giving
the required $\overline{H}$. \epf

The following proposition gives a simple example of S-homotopic morphisms\thinspace:

\bp Let $X$ be a CW-complex and consider the globular CW-complex
$Glob(X)$ (see Proposition~\ref{globe-CW}). Let $x\in X$ and consider
$f_x$ the morphism of globular complexes from $\vI$ to $Glob(X)$ defined
by $f_x(u)=(x,u)$. Take $x$ and $y$ two elements that are in the same
connected component of $X$. Then $f_x$ and $f_y$ are S-homotopic.
\ep

\bpf Let $\gamma$ be a continuous path in $X$ from $x$ to $y$ ($X$
being a CW-complex, its path-connected components coincide with its
connected component).  Then $H(u,t):=(\gamma(t),u)$ satisfies the
condition of Definition~\ref{s-dihomotopic}. \epf

\bd Let $X$ be a globular CW-complex. Then two execution paths
(see Definition~\ref{def_exec_path}) $\gamma_1$ and $\gamma_2$ of $X$
are S-homotopic if and only if the corresponding
morphisms of globular CW-complexes from $\vI$ to $X$ are S-homotopic.
\ed

\bd\label{s-dihomotopic-space} Let $X$ and $Y$ be two
globular CW-complexes.
If there exists
a morphism $f$ from $X$ to $Y$ and a morphism $g$ from $Y$ to $X$
such that $f\circ g\sim_{S} Id_Y$ and $g\circ f\sim_{S} Id_X$, then
$X$ and $Y$ are said \textit{S-homotopic}, or
\textit{S-homotopy equivalent} and $f$ and $g$ are two
inverse S-homotopy equivalences. \ed

Notice that in the latter case, $f\circ g$ and $Id_Y$ coincide on
$Y^0$ and $g\circ f$ and $Id_X$ coincide on $X^0$. Therefore if $f$ is
a S-homotopy equivalence from $X$ to $Y$ then $f$ induces a
bijection between both $0$-skeletons.

Of course, Definition~\ref{s-dihomotopic-space} defines an equivalence
relation.

\bd\label{Sinv} Let $F$ be a functor from $\diCW$ to some category
$\C$. Then $F$ is S-invariant if and only if for any S-homotopy equivalence $s$,
$F(s)$ is an isomorphism in $\C$.  \ed

\bth\label{HOS} Let $S$ be the collection of S-homotopy equivalences of $\diCW$.
There exists a category $\Ho^{S}(\diCW)$ and a functor
$$Q^S:\diCW\longrightarrow \Ho^{S}(\diCW)$$ satisfying the following
conditions\thinspace:
\begin{itemize}
\item For every $s\in S$, $Q^S(s)$ is invertible in $\Ho^{S}(\diCW)$.
\item For every functor $F:\diCW\longrightarrow \C$ such that
for every $s\in S$, $F(s)$ is invertible in $\C$, then there exists a unique
functor $G$ from $\Ho^{S}(\diCW)$ to $\C$ such that $F=G\circ Q^{S}$.
\end{itemize}
\eth

\bpf We mimic the classical proof as presented for instance in
\cite{cube}\thinspace: the main idea consists of using the fact that the
canonical projection from $I^S(X)$ to $X$ is a S-homotopy
equivalence, having as inverse both $\epsilon_0$ and $\epsilon_1$.

Let $\Ho^{S}(\diCW)$ be the category having the same object as $\diCW$
and such that \[\Ho^{S}(\diCW)(X,Y):=[X,Y]_{S}.\] Let
$F:\diCW\longrightarrow \C$ be a functor such that for any $s\in S$,
$Q^{S}(s)$ is invertible in $\C$.  The factorization $F=G\circ Q^{S}$
is obvious on the objects. To complete the proof, it suffices to
verify that for two S-homotopic morphisms $f$ and $g$, then
$F(f)=F(g)$. By definition, there exists $\overline{H}$ from $I^S(X)$
to $Y$ such that $\overline{H}\circ \epsilon_0=f$ and
$\overline{H}\circ \epsilon_1=g$.  Let $pr_1$ be the canonical
projection from $I^S(X)$ to $X$. Then
$pr_1\circ \epsilon_0= pr_1\circ \epsilon_1=Id_X$,
$\epsilon_0\circ pr_1 \sim_{S} Id_{I^S(X)}$ and
$\epsilon_1\circ pr_1 \sim_{S} Id_{I^S(X)}$. Therefore $F(pr_1)$ has
as inverse both $F(\epsilon_0)$ and $F(\epsilon_1)$.  Thus
$F(f)=F(\overline{H})\circ F(\epsilon_0)= F(\overline{H})\circ F(\epsilon_1)=F(g)$.
\epf

\begin{cor}
Let $F$ be a functor from $\diCW$ to some
category $\C$. Then $F$ is S-invariant if and only if there exists a
functor $G$ from $\Ho^{S}(\diCW)$ to $\C$ such that $F=G\circ Q^{S}$.
\end{cor}

\begin{figure}
\begin{center}
\includegraphics[width=7cm]{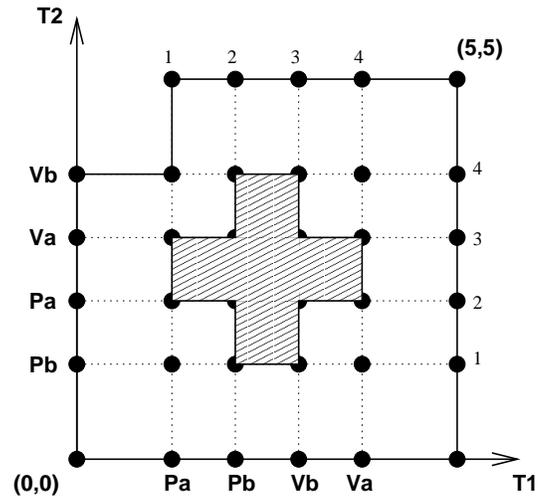}
\end{center}
\caption{Spatial deformation of the HDA of Figure~\ref{progress11}}
\label{progress12}
\end{figure}

\begin{ex}
The HDA of Figure~\ref{progress12} is obtained by contracting the top-left 
$2$-cell of the HDA of Figure~\ref{progress11}. This is a typical example 
of spatial deformation.
\end{ex}

\subsection{T-homotopy equivalence}
\label{Tequiv}

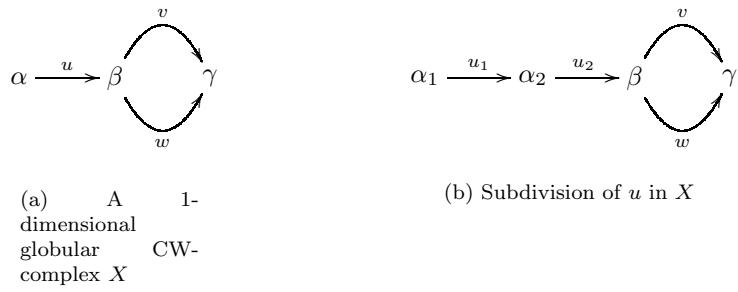
\begin{figure}
\begin{center}
\subfigure[A $1$-dimensional globular CW-complex $X$]{\label{dilate1}
\xymatrix{\alpha\fr{u}& \beta \ar@/^20pt/[r]^v \ar@/_20pt/[r]_w&\gamma}}
\hspace{2cm}
\subfigure[Subdivision  of $u$ in $X$]{
\label{dilate2}
\xymatrix{\alpha_1\fr{u_1}&\alpha_2 \fr{u_2}& \beta \ar@/^20pt/[r]^v \ar@/_20pt/[r]_w&\gamma }}
\end{center}
\caption{Example of temporal deformation}
\label{T1}
\end{figure}

Now we want to treat the case of temporal deformations.
Figure~\ref{T1} is a simple example of temporal deformation of HDA.
An obvious bijective morphism $f$ of globular CW-complexes which sends $u$ on
the ``concatenation'' $u_1 u_2$ (such a morphism is of course not unique)
and which is the identity elsewhere
should be an equivalence. Unfortunately $f$ does not induce a
bijection on the $0$-skeletons because of the point which appears on
the middle of $u$.  So the globular CW-complexes of Figure~\ref{T1}
cannot be S-homotopic. This morphism $f$ induces an homeomorphism
between the underlying topological spaces. The inverse $f^{-1}$ is not
a morphism of globular CW-complexes because the point between $u_1$
and $u_2$ is mapped by $f^{-1}$ on a point belonging to the interior
of the globular cell $u$, which is not an element of the $0$-skeleton.

It is very intuitive to think that morphisms of $\diCW$ inducing
homeomorphisms on the underlying topological spaces do not change
the computer-scientific properties of the corresponding HDA. In
particular, homeomorphisms do not contract directed segments\thinspace:
this is exactly the kind of properties expected for T-invariance
\cite{ConcuToAlgTopo}. Hence the following definition\thinspace:

\bd A morphism $f$ of globular CW-complexes from $X$ to $Y$ is a
\textit{T-homotopy equivalence} if and only if $f$ induces an
homeomorphism on the underlying topological spaces. \ed

\bd\label{Tinv} Let $F$ be a functor from $\diCW$ to some category $\C$. Then
$F$ is T-invariant if and only if for any T-homotopy equivalence  $t$,
$F(t)$ is an isomorphism. \ed

Looking back to Figure~\ref{T1}, one sees that there exists a
T-homotopy equivalence from the left member to the right one, but
not in the reverse direction. So a T-homotopy equivalence is not
necessarily an invertible morphism of $\diCW$.

\bth \label{HOT}
Let $T$ be the collection of T-homotopy equivalences.
There exists a category $\Ho^T(\diCW)$ and a functor
$$Q^T:\diCW\longrightarrow \Ho^T(\diCW)$$ satisfying the following
conditions\thinspace:
\begin{itemize}
\item For every $t\in T$, $Q^T(t)$ is invertible in $\Ho^T(\diCW)$.
\item For every functor $F:\diCW\longrightarrow \C$ such that
for any $t\in T$, $Q^T(t)$ is invertible in $\C$, then there exists a unique
functor $G$ from $\Ho^T(\diCW)$ to $\C$ such that $F=G\circ Q^T$.
\end{itemize}
\eth

\bpf There exists a $\V$-small category $\Ho^T(\diCW)$ satisfying the
universal property of the theorem and constructed as follows\thinspace: the
objects of $\Ho^T(\diCW)$ are those of $\diCW$. The elements of the
$\V$-small set $\Ho^T(\diCW)(X,Y)$ where $X$ and $Y$ are two
$1$-dimensional globular CW-complexes are of the form
\[g_1 f_1^{-1}g_2\dots g_n f_n^{-1} g_{n+1}\]  with $n\geq 1$ where
$g_1,\dots,g_{n+1}$ are morphisms of $\diCW$ and $f_1,\dots,f_n$ are
T-homotopy equivalences and where the notation $f^{-1}$ for $f$
T-homotopy equivalence is a formal inverse of $f$ (see for example
\cite{MR96g:18001a} Proposition 5.2.2 for the construction).

Let us consider the
following commutative diagram
\[
\xymatrix{
X \ar@{->}[r]^(.35){g_{n+1}} & cod(g_{n+1}) & \fl{f_n} dom(f_n) \fr{g_n} & cod(g_n) \\
X \fu{Id}\ar@{->}[r]^(.3){g_{n+1}} & f_n\left(\overline{f_n^{-1}(g_{n+1}(X))}\right)\fu{\subset}
& \ar@{->}[l]_(.45){f_{n}} \fu{\subset}\overline{f_n^{-1}(g_{n+1}(X))} \ar@{->}[r]^(.6){g_{n}} & cod(g_n) \fu{Id}
}
\]
with the notation $cod(h)$ for the codomain of $h$, $dom(h)$ for the domain of $h$,
and for $A$ a subset of some globular CW-complex $Z$,
$$\overline{A}=\bigcup_{x\in A\cap Z^0} \{x\} \cup \bigcup_{x\in
A\backslash Z^0} e_x \subset Z$$ where $e_x$ is the smallest globular cell containing $x$.
We see immediately that
$|\overline{A}|\leq max(2^{\aleph_0},A)$ where $|X|$ means the cardinal of $X$ and where
$\aleph_0$ is the smallest infinite cardinal, i.e. that of the set of non-negative
integers.
Since $f_n$ is an homeomorphism and in particular is bijective, then
\beas
|f_n\left(\overline{f_n^{-1}(g_{n+1}(X))}\right)|&=& |\overline{f_n^{-1}(g_{n+1}(X))}|\\
&\leq & max(2^{\aleph_0},|f_n^{-1}(g_{n+1}(X))|)\\
&=& max(2^{\aleph_0},|g_{n+1}(X)|)\\
&\leq & max(2^{\aleph_0},|X|)
\eeas
This diagram remaining commutative in $\Ho^T(\diCW)$, it shows that
we can suppose $$|cod(g_{n+1})|\leq max(2^{\aleph_0},|X|)$$ and
$$|dom(f_n)|\leq max(2^{\aleph_0},|X|)$$ with an expression like
$g_n f_n^{-1} g_{n+1}$.  By an immediate induction, we see that with a
morphism of the form $g_1 f_1^{-1}g_2\dots g_n f_n^{-1} g_{n+1}$ lying
in $\Ho^T(\diCW)(X,Y)$, we can suppose that all intermediate objects
are of cardinal lower than $max(2^{\aleph_0},|X|)$ which is an
$\U$-small cardinal. Therefore $\Ho^T(\diCW)(X,Y)$ is $\U$-small as well.
\epf

\bp\label{T-invariant} Let $F$ be a functor from $\diCW$ to some category $\C$. Then
$F$ is T-invariant if and only if there exists a functor
$G$ from $\Ho^{T}(\diCW)$ to $\C$ such that
$F=G\circ Q^{T}$. \ep

Let us consider the category $\Ho^{homeo}(\diCW)$ defined as follows\thinspace:
the objects are the globular CW-complexes\thinspace; the set
$\Ho^{homeo}(\diCW)(X,Y)$ is the subset of the set of continuous maps
from $X$ to $Y$ consisting of composites of morphisms of globular
CW-complexes and continuous maps like $f^{-1}$ where $f$ is a
T-homotopy.  There exists a canonical functor $F:\diCW\rightarrow
\Ho^{homeo}(\diCW)$ inverting all T-homotopy equivalences. Therefore there
exists a unique functor $G:\Ho^{T}(\diCW)\rightarrow
\Ho^{homeo}(\diCW)$ such that $G\circ Q=F$.

\begin{question} Why is $G$ an equivalence of categories ? \end{question}

Example~\ref{ex_T} is another typical example of T-homotopy equivalence.

\subsection{Dihomotopy equivalence}
\label{dihomequiv}

Now we want to take into account both spatial and temporal deformations
together.

\bd A morphism of globular CW-complexes is called a dihomotopy
equivalence if it is the composite of S-homotopy equivalence(s) and
T-homotopy equivalence(s).  \ed

\bth \label{HOST}
Let $U$ be the collection of dihomotopy equivalences.
There exists a category $\Ho(\diCW)$ and a functor
$$Q:\diCW\longrightarrow \Ho(\diCW)$$ satisfying the following
conditions\thinspace:
\begin{itemize}
\item For every $u\in U$, $Q(u)$ is invertible in $\Ho(\diCW)$.
\item For every functor $F:\diCW\longrightarrow \C$ such that
for any $u\in U$, $Q(u)$ is invertible in $\C$, then there exists a unique
functor $G$ from $\Ho(\diCW)$ to $\C$ such that $F=G\circ Q$.
\end{itemize}
\eth

\bpf  Let us consider the $\U$-small diagram of categories
\[
\xymatrix{ &  \ar@{->}[ld]_{Q^S}\diCW \ar@{->}[rd]^{Q^T} & \\
\Ho^S(\diCW) & & \Ho^T(\diCW) }
\]
Then the direct limit of this diagram exists in the large category of
$\V$-small categories\thinspace: see \cite{MR96g:18001a} Proposition 5.1.7.  By
reading the construction in the proof of this latter proposition, one
sees that the direct limit is actually a category with $\U$-small
objects and $\U$-small homsets.  \epf

\bp Let $F$ be a functor from $\diCW$ to some
category $\C$. Then $F$ is S-invariant and T-invariant if and only if
there exists a unique functor $G$ from $\Ho(\diCW)$ to $\C$ such that
$F= G\circ Q$. \ep

\bpf Obvious. \epf

\bd The category $\Ho(\diCW)$ is called the category of
dihomotopy types.  \ed

\begin{figure}
\begin{center}
\includegraphics[width=7cm]{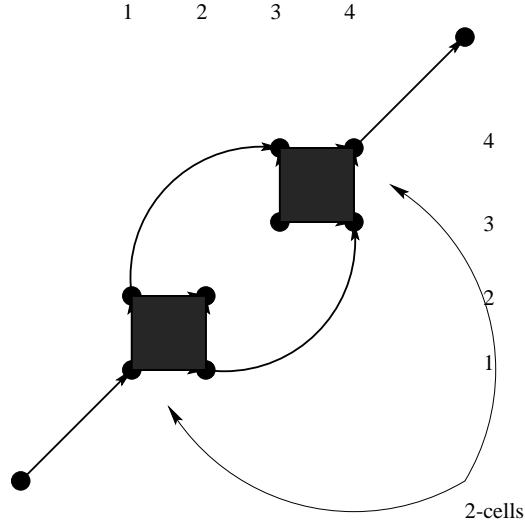}
\end{center}
\caption{HDA dihomotopy equivalent to the HDA of Figure~\ref{progress10}}
\label{progress13}
\end{figure}

\begin{ex} Figure~\ref{progress13} represents the HDA obtained after 
applying several spatial and temporal deformations to the HDA of
Example~\ref{runningexample} (cf. Figure~\ref{progress10}). 
Notice that this HDA has the ``same computer-scientific properties''
as the precubical set from which it was originally built from (see
Example~\ref{runningexample2}). In particular, it has the same execution paths
modulo dihomotopy, and the same set of deadlocks and unreachables.
\end{ex}

\section{Links between homotopy types and dihomotopy types}
\label{embed-ho-diho}

Recall that the category of homotopy types $\Ho(\CW)$ is by definition
the category of CW-complexes with continuous maps up to homotopy, i.e.
if $X$ and $Y$ are two CW-complexes, then $\Ho(\CW)(X,Y):=[X,Y]$. It
is well-known that $\Ho(\CW)$ is the localization of the category
$\CW$ of CW-complexes with respect to the collection of homotopy
equivalences. Theorem~\ref{HOS} can be actually considered as a
generalization of this fact.

\subsection{Path space  between two points}

\label{orthogonal}

Before going any further, we need to define the notion of path
space of a \textit{bipointed local po-space}. Intuitively,
applying this operator to a (global) po-space like $Glob(X)$
(where $X$ is a compactly-generated topological space) bipointed
by $\{\ei,\ef\}$ must give back $X$ up to homotopy.

\bd A \textit{bipointed local po-space} is a triple
$(X,\alpha,\beta)$ where $X$ is a local po-space and $\alpha$ and
$\beta$ are two points of $X$.  A morphism of bipointed local
po-spaces from $(X,\alpha,\beta)$ to $(Y,\alpha,\beta)$ is a
morphism of po-spaces $f$ from $X$ to $Y$ such that
$f(\alpha)=\alpha$ and $f(\beta)=\beta$. The corresponding
category is denoted by $\lpohauspp$. \ed

Notice that $Glob(-)$ can be seen as a functor from $\haus$ to
$\lpohauspp$ or from $\top$ to $\potoppp$ (the category of bipointed
topological spaces with a non-necessarily closed partial ordering) by
bipointing $Glob(X)$ by the elements $\ei$ and $\ef$.

\bp\label{dlim} The functor $Glob(-)$ from $\top$ to $\potoppp$ commutes
with direct limits. \ep

\bpf Let $(X_i)_{i\in I}$ be a family of topological spaces.
Then
\[Glob\left(\bigsqcup_{i\in I} X_i\right)=
\bigsqcup_{i\in I} (X_i\p[0,1])
\left/
\left\{\begin{array}{c}
(z,0)=(z',0) \hbox{ for }z,z'\in \bigsqcup_{i\in I} X_i\\
(z,1)=(z',1) \hbox{ for }z,z'\in \bigsqcup_{i\in I} X_i
\end{array}\right\}\right.\]

Note that for all
$x\in Glob\left(\bigsqcup_{i\in I} X_i\right)\backslash
\{\alpha,\beta\}$, there exists a unique $i_x\in I$ such that
$x\in Glob(X_{i_x})$.
Let $(T,\alpha,\beta)$ be a bipointed topological space and
for all $i\in I$, let
$\phi_i:Glob(X_i)\longrightarrow (T,\alpha,\beta)$
be a morphism in $\potoppp$. Let $\phi$ be the set map from
$Glob\left(\bigsqcup_{i\in I} X_i\right)$ to
$T$ defined by
$\phi(\alpha)=\alpha$, $\phi(\beta)=\beta$, and
$\phi(x)=\phi_{i_x}(x)$ (for $x \neq \alpha$ and $x \neq \beta$).
Take $(x,t), (x,t') \in Glob\left(\bigsqcup_{i\in I} X_i\right)$
such that $(x,t)\leq (x',t')$
We have three possibilities:
\begin{itemize}
\item
$(x,t)=\alpha$ and $\phi(x,t)=\alpha\leq \phi(x',t')$,
\item
$(x',t')=\beta$ and $\phi(x,t)\leq \beta = \phi(x',t')$,
\item
$\alpha<(x,t)\leq (x',t')<\beta$.
\end{itemize}
In the latter case, $x=x'$ and
therefore there exists $i_0\in I$ such that $(x,t)$ and $(x',t')$
belong to $Glob(X_{i_0})$.  Then
$\phi(x,t)=\phi_{i_0}(x,t)\leq \phi_{i_0}(x',t')=\phi(x',t')$.  The
set map $\phi$ is well-defined and continuous because it is the
quotient in $\top$ of the direct sum $\bigsqcup_{i\in I}\phi_i$ by the
identifications $(z,0)=(z',0) \hbox{ for }z,z'\in \bigsqcup_{i\in I}
X_i$ and $(z,1)=(z',1) \hbox{ for }z,z'\in \bigsqcup_{i\in I} X_i$.
Therefore $Glob\left(\bigsqcup_{i\in I} X_i\right)$ is the
direct sum of the $Glob(X_i)$ for $i$ running over $I$
in $\potoppp$. So  the functor $Glob(-)$ preserves the direct
sums.

Let $f$ and $g$ be two continuous maps from $X$ to $Y$. Let
$Z=Y\left/ \{ f(x) \equiv g(x) \mid x\in X\}\right.$ be the
coequalizer of $(f,g)$ in $\top$. Then there exists a surjection
{\small
\[
Glob(Y) \left/\left\{
Glob(f)(x,t) \equiv Glob(g)(x,t)
\right\}\right.
\twoheadrightarrow
\left( \left(Y\left/ \{f(x)=g(x)\}\right.\right) \p [0,1]\right)
\left/ \left\{\begin{array}{c}
(z,0)=(z',0) \\
(z,1)=(z',1)
\end{array}\right\}\right.
\]}
which is clearly an homeomorphism.
Let $(T,\alpha,\beta)$ be a bipointed
local po-space and let $h$ be a morphism in $\potoppp$ from
$Glob(Y)$ to $T$ such that
$h\circ Glob(f)=h\circ Glob(g)$.
Then $h$ factorizes through $Glob(Z)$ because this latter
is the coequalizer of $(Glob(f),Glob(g))$
in $\top$. It is
easily checked that $h$ is a non-decreasing map and therefore a
morphism. So $Glob(-)$ preserves the coequalizers. This entails the result
by Proposition 2.9.2 of \cite{MR96g:18001a}.
\epf

\bp\label{radjoint} The functor $Glob(-)$ from $\top$ to $\potoppp$
has a right adjoint that will be denoted by $(-)^{\bot}$.  \ep

\bpf The category $\top$ has a generator\thinspace: the one-point space;
it is cocomplete and well-copowered. The result follows from the
Special Adjoint Functor theorem \cite{MR1712872}. \epf

If $X$ and $Y$ are two topological spaces, the topological space
$Cop(X,Y)$ will be by definition the set $\top(X,Y)$ of continuous
maps from $X$ to $Y$ endowed with the compact-open topology\thinspace: a basis
for the open sets consists of the sets $N(C,U)$ where $C$ is any
compact subset of $X$, U any open subset of $Y$ and
\(N(C,U):=\{f\in\top(X,Y),f(C)\subset U\}\).  A topological space is
compactly generated when its topology coincides with the weak topology
determined by its compact subspaces. If $X$ is a Hausdorff topological
space, its Kelleyfication $k(X):=\limind_{K\subset X} K$ is a
Hausdorff compactly generated topological space where $K$ runs over
the set of compact subspaces of $X$. The Kelleyfication functor $k$ is
the right adjoint of the forgetful functor from compactly generated
topological spaces to Hausdorff topological spaces. So if $X$ is
compactly generated and if $Y$ is Hausdorff, then a set map
$f:X\rightarrow Y$ is continuous if and only if $f:X\rightarrow k(Y)$
is continuous. The topological space $Cop(X,Y)$ needs not be compactly
generated. Let $\thaus(X,Y):=k(Cop(X,Y))$.  Every locally compact
Hausdorff topological space is compactly generated. In particular,
every CW-complex and every globular CW-complex is compactly generated.
The main property of the compact-open topology is the following one\thinspace:
If $X$, $Y$ and $Z$ are compactly generated, then one has a natural
bijection of sets
\begin{equation}\label{cartesian}
\top(X\p Y,Z)\iso \top (X,\thaus(Y,Z))
\end{equation}
induced by $f\mapsto (x\mapsto f(x,-))$ from the left to the right
member and by $g\mapsto ((x,y)\mapsto g(x)(y))$ in the opposite
direction. As a matter of fact, the isomorphism (\ref{cartesian}) as
topological spaces holds
as soon as $Y$ is locally compact Hausdorff.

\bp
\label{prop43}
If $(X,\alpha,\beta)$ is a bipointed local po-space such that
$X$ is compactly generated, then $(X,\alpha,\beta)^{\bot}$ is homeomorphic to
the set of non-decreasing maps $\gamma$ from $[0,1]$ to $X$ such that
$\gamma(0)=\alpha$ and $\gamma(1)=\beta$, endowed with the
Kelleyfication of the compact-open topology. \ep

\bpf Since $[0,1]$ is compact,
\[\top(Y\p [0,1],X)\iso \top (Y,\thaus([0,1],X)).\]
This isomorphism specializes to
\[
\left\{f:Y\p [0,1]\longrightarrow X,
\begin{array}{c}
f(y,0)=\alpha \hbox{ for all }y\in Y\\
f(y,1)=\beta \hbox{ for all }y\in Y\\
f(y,-)\hbox{ dipath of }X
\end{array}
\right\}
\iso
\left\{g:Y\longrightarrow \thaus_{\alpha,\beta}([0,1],X)\right\}
\]
where $\thaus_{\alpha,\beta}([0,1],X)$ is the set of non-decreasing
continuous maps $\gamma$ from $[0,1]$ to $X$ such that
$\gamma(0)=\alpha$ and $\gamma(1)=\beta$.
The first member is in natural bijection with the morphisms
of bipointed po-spaces from $Glob(Y)$ to $X$, hence the result.
\epf

\bd\label{ortho}
\label{def42}
For $(X,\alpha,\beta)\in\potoppp$ with $X$ compactly generated,
then the topological space
\[\P(X,\alpha,\beta):=(X,\alpha,\beta)^{\bot}\backslash\{\alpha\}\]
where $\alpha$ is the constant path $\alpha$, is called the
\textit{path space} of $(X,\alpha,\beta)$, or the \textit{path space }
of $X$ from $\alpha$ to $\beta$. Notice that
$\alpha\in (X,\alpha,\beta)^{\bot}$ if and only if $\alpha=\beta$. \ed

If $\P Glob(X)$ is only equipped with the compact-open topology, it will
be rather denoted by $\P^{co} Glob(X)$. In other terms,
$\P Glob(X)=k(\P^{co} Glob(X))$.

The canonical map $i$ from $X$ to $\P Glob(X)$
maps any $x\in X$ to the dimap $t\mapsto (x,t)$ of
$\P Glob(X)$. Now,

\bth\label{homotopy}
\label{thm44}
For any compactly generated topological space $X$, the canonical map
from $X$ to $\P Glob(X)$ is an homotopy
equivalence. \eth

\bpf Let $\phi \in \P Glob(X)$.
By definition, $\phi$ is a non-decreasing continuous path from
$\phi(0)=\ei$ to $\phi(1)=\ef$. Let $pr_2$ be the canonical
projection of
$Glob(X)$ onto $[0,1]$. Since $]0,1[$ is open and connected, and
$pr_2$ and $\phi$ are continuous, $(pr_2\circ
\phi)^{-1}(]0,1[)$ is open and connected. Thus we can set $(pr_2\circ
\phi)^{-1}(]0,1[)=]t_{\phi}^-,t_{\phi}^+[$.  Due to the peculiar ordering we
have on $Glob(X)$, $\phi$ being non-decreasing implies that there exists a unique
$\underline{x}(\phi)\in X$ such that for $t\in
]t_{\phi}^-,t_{\phi}^+[$,
$\phi(t)=(\underline{x}(\phi),pr_2\circ\phi(t))$ (i.e. its first component
is constant on $]t_{\phi}^-,t_{\phi}^+[$).

Let $\phi_0\in \P^{co} Glob(X)$ and let $U$ be an open of X
containing $\underline{x}(\phi_0)$. Let $K_{\phi_0}$ be a compact
subset of $]t_{\phi_0}^-,t_{\phi_0}^+[$. Then $\phi_0 \in
N(K_{\phi_0},U\p ]0,1[)$  and for every $\phi\in N(K_{\phi_0},U\p
]0,1[)$, $\underline{x}(\phi)\in U$.  Therefore the map
$\underline{x}$ from $\P^{co} Glob(X)$ to $X$ is continuous. Therefore
the map still denoted by $\underline{x}$ from $\P Glob(X)$ to $k(X)=X$
is continuous.

One has $\underline{x}\circ i=Id_X$ and for all $\phi\in \P Glob(X)$,
$i\circ \underline{x}(\phi)$ is the dimap
$t\mapsto (\underline{x}(\phi),t)$. Let
\[H(\phi,u)(t)=(\underline{x}(\phi),ut+(1-u)pr_2\circ \phi(t))\]
Then $H$ yields a set map from $\P Glob(X) \p \I$ to
$\P Glob(X)$ with $H(\phi,0)=\phi$ and $H(\phi,1)=i\circ
\underline{x}(\phi)$.  So it suffices to check the continuity of $H$
to complete the proof.

Consider the set map $H'$ from
$\P^{co} Glob(X)\p \I\p \I$ to
$\P^{co} Glob(X)$ defined by
\[H'(\phi,u,t)=(\underline{x}(\phi),ut+(1-u)pr_2 \circ \phi(t))\]
Let $C$ be a compact subset of $\I$ and $U$ be an open subset of
$\I$ such that $pr_2\circ \phi_0\in N(C,U)$ for some
$\phi_0\in \P^{co} Glob(X)$. Then for any $\phi\in
N(C,X\p U)$, $pr_2\circ \phi(C)\subset U$. Therefore
the set map $pr_2:\P^{co} Glob(X)\longrightarrow \thaus(\I,\I)$
defined by $pr_2(\phi)=pr_2 \circ \phi$ is
continuous, and the set map $H'$ is continuous as well. Since the
Kelleyfication functor is a right adjoint, it then commutes with
products. So $k(H')$ is a continuous map from
$\P Glob(X)\p \I\p \I$ to $\P Glob(X)$.

Since $H$ is the image of $k(H')$ by the canonical isomorphism

{\small
\[\top(\P Glob(X)\p
\I\p \I,\P Glob(X))
\longrightarrow \top (\P Glob(X)\p
\I,\thaus(\I,\P Glob(X)))
\] }

$H$ is continuous as well.
\epf

\subsection{Homotopy and dihomotopy types}

We have now the necessary tools in hand to compare homotopy types and
dihomotopy types.

\bth\label{th} Let $X$ and $Y$ be two compactly generated topological spaces.
Let $f$ be a morphism of globular complexes from $Glob(X)$ to $Glob(Y)$.
Then there exists a unique continuous map $f^S$ from $X$ to $Y$ up to
homotopy such that $f$ is S-homotopic to $Glob(f^S)$. \eth

\bpf Let $f_0$ and $f_1$ be two continuous maps from $X$ to $Y$ such
that $Glob(f_0)$ and $Glob(f_1)$ are S-homotopic to $f$. Let $H$ from
$Glob(X)\p \I$ to $Glob(Y)$ be a S-homotopy from $Glob(f_0)$ to $Glob(f_1)$
with $H_t:=H(-,t)$, $H_0=Glob(f_0)$ and $H_1=Glob(f_1)$. Consider the set
map $h$ from $X\p \I$ to $Y$ defined by
$h(x,t)=(\underline{x}\circ \P(H_t) \circ i)(x)$ with the notations
of Theorem~\ref{homotopy}. Then
\beas
h(x,0)&=&(\underline{x}\circ \P(h_0) \circ i)(x)\\
&=&(\underline{x}\circ \P Glob(f_0))\left(u\mapsto (x,u)\right)\\
&=&\underline{x}\left(u\mapsto (f_0(x),u)\right)\\
&=& f_0(x) \eeas and in the same manner one gets $h(x,1)=f_1(x)$.
So it suffices to prove the continuity of $h$ to prove the
uniqueness of $f^S$ up to homotopy. We have already proved in
Theorem~\ref{homotopy} the continuity of $i$ and
$\underline{x}$.  Therefore it suffices to prove the continuity
of the set map $(\gamma,t)\mapsto \P(H_t)(\gamma)=H_t\circ
\gamma$ from $\P Glob(X)$ to $\P Glob(Y)$. This latter map is the
composite of
\[
\xymatrix{ {\P Glob(X) \p \I}\fd{(\gamma,t)\mapsto (\gamma,H,t)}\\
{\P Glob(X) \p \thaus(Glob(X)\p \I,Glob(Y))\p \I}\fd{(\gamma,H,t)\mapsto (\gamma,H_t)}\\
{\P Glob(X) \p \thaus(Glob(X),Glob(Y))}\fd{(\gamma,g)\mapsto g\circ \gamma}\\
{\P Glob(Y)}}
\]
The last map $(\gamma,g)\mapsto g\circ \gamma$ is the image of
the identity map of $\thaus(Glob(X),Glob(Y))$ by
{\footnotesize
\[
\xymatrix{{\top(\thaus(Glob(X),Glob(Y)),\thaus(Glob(X),Glob(Y)))}\fd{\iso}\\
{\top(Glob(X)\p \thaus(Glob(X),Glob(Y)),Glob(Y))}\fd{}\\
{\top(\vI\p \thaus(\vI,Glob(X)) \p \thaus(Glob(X),Glob(Y)),Glob(Y))}\fd{\iso}\\
{\top(\thaus({\vI},Glob(X))\p \thaus(Glob(X),Glob(Y)),\thaus({\vI},Glob(Y)))}}
\]}
and therefore is continuous. At last the set map $(H,t)\mapsto H_t$
is the image of the identity map of $\thaus(Glob(X)\p \I,Glob(Y))$ by
{\footnotesize
\[
\xymatrix{{\top(\thaus(Glob(X)\p \I,Glob(Y)),\thaus(Glob(X)\p \I,Glob(Y)))}\fd{}\\
{\top(Glob(X)\p \I \p \thaus(Glob(X)\p \I,Glob(Y)),Glob(Y))}\fd{\iso}\\
{\top(\I \p \thaus(Glob(X)\p \I,Glob(Y)),\thaus(Glob(X),Glob(Y)))}
}
\]}
and therefore is also continuous. So $h$ is an homotopy between
$f_0$ and $f_1$.

Now set $f^S:=\underline{x}\circ \P(f) \circ i$ from $X$ to $Y$. With
the proof of Theorem~\ref{homotopy}, we see immediately that $f^S$ is
continuous.  It remains to prove that $Glob(f^S)$ is S-homotopic to $f$.
We have already seen in the proof of
Theorem~\ref{homotopy} that for $\phi\in \P Glob(X)$,
\begin{equation}\label{decomposition}
\phi(t)=(\underline{x}(\phi),pr_2\circ\phi(t))
\end{equation}
 for $t\in ]t_{\phi}^-,t_{\phi}^+[$. For $t\in [0,t_{\phi}^-]$ (resp.
 $t\in [t_{\phi}^+,1]$), one has by definition $pr_2\circ\phi(t)=0$
 (resp. $pr_2\circ\phi(t)=1$) and therefore Equality~\ref{decomposition} is
 still true for any $t\in \vI$. So consider the path $\phi_x:t\mapsto (x,t)$
of $\P Glob(X)$. Then $f\circ\phi_x$ is an element of $\P Glob(Y)$ and we have
$f\circ\phi_x=(\underline{x}(f\circ\phi_x),pr_2\circ f\circ\phi_x (t))$.
But $\underline{x}(f\circ\phi_x)=f^S(x)$. Therefore $f=(f^S,pr_2\circ f)$. So
$f$ is S-homotopic to $Glob(f^S)$ with the S-homotopy $H$ from $Glob(X)\p \I$
to $Glob(Y)$ defined by $H((x,t),u)=(f^S(x),ut+(1-u)pr_2\circ f(x,t))$.
\epf

\begin{cor}\label{bij-hom-dihom} Let $X$ and $Y$ be two compactly generated
topological spaces.  The functor $Glob(-)$ induces a bijection of sets
$[X,Y]\iso [Glob(X),Glob(Y)]_{S}$.
\end{cor}

We arrive at

\bth\label{embed} The mapping $X\mapsto Glob(X)$ induces an embedding
$$\Ho(\CW)\hookrightarrow \Ho(\diCW).$$  \eth

\bpf This is a consequence of Proposition~\ref{globe-CW}, Proposition~\ref{globe-CW2}
and Theorem~\ref{th}. \epf

See the consequences of this important theorem in \cite{Ditype} where a research
program to investigate dihomotopy types is exposed.

\subsection{Towards a Whitehead theorem}

Now we want to weaken the notion of S-homotopy equivalence.

\bd Let $f$ be a morphism of globular CW-complexes from $X$ to $Y$. Then
$f$ is a weak S-homotopy equivalence if the following conditions are fulfilled\thinspace:
\begin{enumerate}
\item the map $f$ induces a set bijection between the $0$-skeleton of $X$ and
the $0$-skeleton of $Y$.
\item for $\alpha,\beta\in X^0$, $f$ induces a weak homotopy equivalence
from $\P (X,\alpha,\beta)$ to $\P (Y,f(\alpha),\linebreak[4]f(\beta))$.
\end{enumerate}
\ed

\bp\label{strong-implies-weak} Let $f$ be a morphism of globular
CW-complexes from $X$ to $Y$. If $f$ is a S-homotopy from $X$ to
$Y$, then $f$ is a weak S-homotopy equivalence. \ep

\bpf Let $g$ be a S-homotopy from $Y$ to $X$ such that $f\circ
g\sim_{di} Id_Y$ and $g\circ f\sim_{di} Id_X$. Then $f\circ g$
and $Id_Y$ (resp. $g\circ f$ and $Id_X$) coincide on $Y^0$ (resp.
$X^0$). Therefore $f$ induces a bijection of sets from the
$0$-skeleton $X^0$ to the $0$-skeleton $Y^0$ with inverse the
restriction of $g$ to $Y^0$. Let $\alpha$ and $\beta$ be two
elements of $X^0$. Then $f$ (resp. $g$) induces a continuous map
$f_*$ from $\P(X,\alpha,\beta)$ (resp. $g_*$ from
$\P(Y,f(\alpha),f(\beta))$) to $\P(Y,f(\alpha),f(\beta))$ (resp.
$\P(X,\alpha,\beta)$). Let $H$ be a continuous map from $Y\p \I$
to $Y$ which is a S-homotopy from $f\circ g$ to $Id_Y$. Let
$H_u=H(-,u)$. By hypothesis, this is a morphism of globular
CW-complexes from $Y$ to itself which induces the identity map on
$Y^0$. Let $h(\gamma,u):=H_u\circ \gamma$. Then
$h(\gamma,u)(0)=H_u(\gamma(0))=H_u(f(\alpha))=f(\alpha)$ and
$h(\gamma,u)(1)=H_u(\gamma(1))=H_u(f(\beta))=f(\beta)$. Moreover
$h(\gamma,u)$ is non-decreasing and continuous because it is the
composite of two functions which are non-decreasing and
continuous as well.  Therefore $h$ is a set map from
$\P(Y,f(\alpha),f(\beta))\p \I$ to $\P(Y,f(\alpha),f(\beta))$. We
have already proved the continuity of similar maps (as in
Theorem~\ref{th}). Therefore $f_*\circ g_* \sim
Id_{\P(Y,f(\alpha),f(\beta))}$. Similarly, we can prove that
$g_*\circ f_* \sim Id_{\P(X,\alpha,\beta)}$. Therefore $f$ is a
weak S-homotopy equivalence. \epf

The converse of Proposition~\ref{strong-implies-weak} gives rise to the following

\begin{conj}\label{strong-eq-weak}\footnote{This conjecture has been actually 
solved later on \cite{Wh-Gau}.}

 Let $f$ be a morphism of globular CW-complexes
from $X$ to $Y$. Then the following assumptions are equivalent\thinspace:
\begin{enumerate}
\item $f$ is a weak S-homotopy equivalence.
\item $f$ is a S-homotopy equivalence.
\end{enumerate}
\end{conj}

In the case of globes, one has\thinspace:

\bp Let $f$ be a morphism of globular CW-complexes from $Glob(X)$ to
$Glob(Y)$ where $X$ and $Y$ are two connected CW-complexes. If $f$ is a
weak S-homotopy equivalence, then there exists a morphism of globular
CW-complexes $g$ from $Glob(Y)$ to $Glob(X)$ such that $g\circ f$ is S-homotopic
to the identity of $Glob(X)$ and $f\circ g$ S-homotopic to the identity of
$Glob(Y)$. \ep

\bpf The composite $\underline{x}\circ \P(f) \circ i$
\[\xymatrix{X\fr{}& {\P Glob(X)} \ar@{->}[rr]^{\P(f)} && {\P Glob(Y)} \fr{}& Y}\]
is a homotopy  equivalence of CW-complexes because $\P(f)$ is
an homotopy equivalence by hypothesis and because of
Theorem~\ref{homotopy}. Therefore
$\underline{x}\circ \P(f) \circ i$ has an inverse $g$ up to
homotopy from $Y$ to $X$. By Corollary~\ref{bij-hom-dihom},
$Glob(\underline{x}\circ \P(f) \circ i)\circ Glob(g)$ and
$Glob(g)\circ Glob(\underline{x}\circ \P(f) \circ i)$ are S-homotopic to the identity
(resp. of $Glob(Y)$ and $Glob(X)$). Again by Corollary~\ref{bij-hom-dihom},
$Glob(\underline{x}\circ \P(f) \circ i)$ and $f$ are S-homotopic. Therefore
$$Glob(g)\circ f \sim_{S} Glob(g)\circ Glob(\underline{x}\circ \P(f) \circ i)\sim_{S} Id$$ and
$$f\circ Glob(g)\sim_{S} Glob(\underline{x}\circ \P(f)  \circ i)\circ Glob(g)\sim_{S} Id.$$
\epf

\section{Why non-contracting maps ?}
\label{why}

We would like to explain here why one imposes the morphisms of
globular CW-complexes to be non-contracting in
Definition~\ref{morphism-glCW}, why in Definition~\ref{ortho} the
constant dipath is removed from $\P(X,\alpha,\beta)$ if
$\alpha=\beta$. As a matter of fact, there are a lot of technical
reasons to do that which will be clearer in the future developments.
This section focuses on a very striking one.

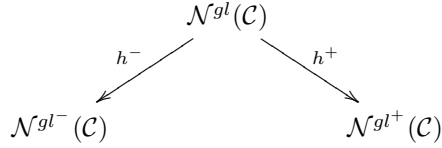
\begin{figure}
\[\xymatrix{& {\mathcal{N}^{gl}(\C)}\ar@{->}[ld]_{h^-}\ar@{->}[rd]^{h^+}&\\
{\mathcal{N}^{gl^-}(\C)}&&{\mathcal{N}^{gl^+}(\C)}}\]
\caption{The fundamental diagram}
\label{fundamental}
\end{figure}

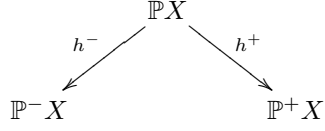
\begin{figure}
\[\xymatrix{& {\P X}\ar@{->}[ld]_{h^-}\ar@{->}[rd]^{h^+}&\\
{\P^- X}&&{\P^{+}X}}\]
\caption{The fundamental diagram for a globular CW-complex $X$}
\label{fundamental2}
\end{figure}

The fundamental algebraic structure which has emerged from the
$\omega$-categorical approach \cite{Gau,sglob,math.AT/0103011} is the
diagram of Figure~\ref{fundamental} where $\C$ is an
$\omega$-category. The analogue in the globular CW-complex framework
is the diagram of Figure~\ref{fundamental2} where $\P X$ is the space
of dipaths between two elements of the $0$-skeleton of $X$, and $\P^-
X$ (resp. $\P^+ X$) is the space of germs of dipaths starting from
(resp. ending at) a point of the $0$-skeleton of $X$.

Let us suppose just for this section  that the \textit{path space} of a globular
CW-complex $X$ is defined as follows\thinspace:
\[ \P X = \bigsqcup_{(\alpha,\beta)\in X^0\p X^0} (X,\alpha,\beta)^\bot \]
and that the \textit{semi-path spaces} of a globular CW-complex $X$ are
defined as follows\thinspace:
 \beas
&& \P^- X = \bigsqcup_{\alpha\in X^0} (X,\alpha)^{\bot^-} \\
&& \P^+ X = \bigsqcup_{\alpha\in X^0} (X,\alpha)^{\bot^+} \eeas
where $(X,\alpha)^{\bot^-}$ (resp. $(X,\alpha)^{\bot^+}$) is the set
of dipaths of $X$ starting from $\alpha$ (ending at $\alpha$), and
endowed with the compact-open topology. Then the maps $h^-$ and $h^+$ of
Figure~\ref{fundamental2} are obviously defined. However

\bp (Remark due to Stefan Sokolowski) The topological spaces
$\P^- X$ and $\P^+ X$ are homotopy equivalent to the discrete set
$X^0$ (the $0$-skeleton of $X$) ! \ep

\bpf Let us make the proof for $\P^- X$. The canonical map $u:X^0
\hookrightarrow \P^- X$ sends an $\alpha\in X^0$ on the corresponding
constant dipaths of $\P^- X$. The map $u$ is necessarily continuous
since $X^0$ is discrete. In the other direction, let us consider the
set map $v:\P^- X \rightarrow X^0$ defined by $v(\gamma)=\gamma(0)$\thinspace:
such an evaluation map is necessarily continuous as soon as $\P^- X$
is endowed with the compact-open topology. Then $v\circ u=Id_{X^0}$
and $u\circ v$ is homotopic to $Id_{\P^- X}$ by the homotopy
\[H: \P^- X \p \I \rightarrow \P^- X\]
defined by $H(\gamma,u)(t):=\gamma(tu)$. The map $H$ is the image of
the identity of $\thaus(\vI,X)$ by
\[
\xymatrix{
\top (\thaus(\vI,X), \thaus(\vI,X)) \fd{\iso}\\
\top(\vI\p \thaus(\vI,X),X)\fd{\phi}\\
\top(\vI\p \thaus(\vI,X)\p\I,X)\fd{\iso}\\
\top(\thaus(\vI,X)\p\I,\thaus(\vI,X))
}
\]
where $\phi$ is induced by the mapping $(t,u)\mapsto tu$ from $\vI\p \I$ to
$\vI$ and therefore $H$ is continuous.
\epf

Therefore $\P^- X$ and $\P^+ X$ defined as above contain no relevant
information !  This fact is exactly the analogue of \cite{Gau}
Proposition 4.2 which states that the cubical nerve of an
$\omega$-category has a trivial simplicial homology with respect to
$\de^-$ and $\de^+$ and which led  to introducing $\omega
Cat(I^*,\C)^-$ and $\omega Cat(I^*,\C)^+$.

So the correct definition of $\P^- X$ and $\P^+ X$ is\thinspace:

 \beas
&& \P^- X = \bigsqcup_{\alpha\in X^0} (X,\alpha)^{\bot^-}\backslash\{\alpha\} \\
&& \P^+ X = \bigsqcup_{\alpha\in X^0} (X,\alpha)^{\bot^+}\backslash\{\alpha\} \eeas

Now the maps $\P X\rightarrow \P^- X$ and $\P X\rightarrow \P^+ X$ do not
exist anymore ! To recover these important maps, it is necessary to set\thinspace:
\[ \P X = \bigsqcup_{(\alpha,\beta)\in X^0\p X^0} \P (X,\alpha,\beta) \]

Then the only way to make the mapping $\P$ (and also $\P^-$ and
$\P^+$) a functor from the category $\diCW$ of globular CW-complexes
to that of compactly-generated topological spaces is to impose to
morphisms in $\diCW$ to be non-contracting as explained in
Definition~\ref{morphism-glCW}.

\section{Concluding remarks and some open questions}

We have constructed a category of dihomotopy types whose isomorphism
classes of objects represent exactly higher dimensional automata
modulo deformations leaving invariant computer-scientific properties
as presence or not of deadlock or everything related.  This
construction provides a rigorous definition of S-deformations
(Definition~\ref{Sinv}) and T-deformations (Definition~\ref{Tinv}) of
HDA. Using the definitions of \cite{Ditype}, it is trivial to prove
the S-invariance of all functors like $H_*^{gl}$, $H_*^{gl^\pm}$, etc...

\begin{question} Proving the T-invariance of both semi-globular homology
theories $H_*^{gl^\pm}$. 
\end{question}

\begin{question} Same question for the biglobular homology defined in \cite{sglob}.
\end{question}

By analogy with the situation in usual algebraic topology\thinspace:

\begin{question} Defining a notion of weak dihomotopy equivalence on the category
of local po-spaces\thinspace; Proving that the localization of the category of local po-spaces
with respect to this collection of morphisms exists and that it is isomorphic
to the category of dihomotopy types. \end{question}

The realization functor from a (quite large) subcategory of precubical sets
(the ``non-self-linked'' ones) to the category of local po-spaces constructed
in \cite{LFEGMRAlgebraic} and the realization functor constructed in
Section~\ref{appA} must be compared, which leads to the following

\begin{question} Proving that for a non-self linked precubical set,
both realization functors give the same local po-spaces up to
dihomotopy\thinspace; so a notion of dihomotopic local po-spaces is needed.
\end{question}

In sequential computation theory, there are nice algebraic topological
results which should relate to the present theory.

A mono\"{\i}d $(M,1,.)$ (1 is the neutral element for the monoidal operation .)
is finitely presented (respectively finitely presented
by a rewriting system) if there exists a finite set of symbols $S$ and
a finite set of relations $R$ on $S^*$, the free monoid on $S$
(respectively of directed rewrite rules of the
form $u$ rewrites into $v$, $u$, $v$ in $S^*$), such that $M \cong S^*/R$ (respectively
such that $M$ is isomorphic to $S^*$ quotiented by the congruence generated
by $R$). We say that a mono\"{\i}d $(M,1,.)$ has a decidable word
problem if there exists an algorithm which can decide in a finite time
whether two words of $M$ are equal.

We say that a rewriting system is canonical if given any word $w$ in $S^*$,
any sequence of reductions (i.e. sequence of applications of any of the rewrite
rules) is finite, and if all sequences of reductions (since many different
reductions can apply at the same time on the same word) end at the same word.
This last word is called the canonical form of $w$ and is used to decide equality
of two representatives of words in $M$.

Squier's theorem shows that all finitely presented mono\"{\i}ds
with a decidable word problem cannot be presented by a finite canonical
rewriting system. In fact, it tells us a lot more: if $M$ can be presented
by a finite canonical rewriting system, then the third homology group
of $M$ is of finite type. In fact, we can even prove \cite{YKComplete,CCSFOYKFiniteness}
that all homology groups are of finite type.
One of the resolutions \cite{JRJGRewriting} constructs a cubical set on which
the mono\"{\i}d acts freely on the left. This means that all orbits of points
under the action of $M$ are trajectories on a cubical set; then, the homology
of the mono\"{\i}d is the homology of the quotient space of the geometric realization
of the cubical set constructed under the relation, being on the same orbit.

It is well known that in the case of mono\"{\i}ds, having finite type homology
groups does not imply the existence of a finite canonical rewrite system presenting
the mono\"{\i}d, whereas in the case of groups, more things are known. This
``defect'' is partly due to the fact that we forget in the classical homological
construction the ``direction of time'' represented by the (non-commutative and
non-invertible in
general) mono\"{\i}dal operation. It is also probably due to the fact that
the real criterion would have to be homotopical and not just homological.
These ideas are already reflected in
\cite{CCSFOYKFiniteness}
and in the more recent one
\cite{YLNew}. In fact, in the case of concurrency theory,
we only deal with ``partially commutative'' mono\"{\i}ds i.e. free mono\"{\i}ds
modulo the commutativity of certains pair of actions, or of certains $n$-uples
of actions. This is only a particular case, which is more general in the
case of rewriting systems.
As a matter of fact, in \cite{YLNew}, the
relations of the mono\"{\i}d are presented by 2-cells in a 2-category
(i.e. homotopies between [di-]paths!).

In modern terms, the homology of a mono\"{\i}d is the homology of its classifying
space. The classifying space is just the geometric realization of the nerve of
the mono\"{\i}d (considered as a one-object category). Our point of view is that
the globular homology of the mono\"{\i}d (i.e. the homology of the globular
nerve of the mono\"{\i}d) is an ``invariant'' of the mono\"{\i}d, i.e. will
not change when we change the presentation (change the set of symbols and relations).
The problem is that this globular homology is far too big (in general of infinite
rank) for our purposes. For it to be reasonable, it would have to be somehow
invariant under subdivisions. This is currently under investigation.

\section*{Acknowledgments}

The authors thank Gunnar Carlsson and Rick Jardine for the invitation
to the conference ``Algebraic Topological Methods in Computer
Science'' held in Stanford in August 2001.  We also wish to thank the
anonymous referees and Martin Raussen for helpful comments.

\end{document}